\documentclass[11pt]{article}
\usepackage[a4paper,margin=2.4cm]{geometry}
\usepackage{amsmath,amssymb,amsthm,booktabs,array,graphicx,multirow}
\usepackage{placeins}
\usepackage[numbers,sort&compress]{natbib}
\usepackage{xcolor}
\usepackage{hyperref}
\usepackage[T1]{fontenc}
\definecolor{linkblue}{rgb}{0,0.2,0.45}
\hypersetup{colorlinks=true,linkcolor=linkblue,citecolor=linkblue,urlcolor=linkblue}
\newcommand{\doi}[1]{\url{https://doi.org/#1}}

\graphicspath{{figures/}}
\newtheorem{proposition}{Proposition}

\newtheorem{remark}{Remark}
\newcommand{\avg}[1]{\{\!\!\{#1\}\!\!\}}
\newcommand{\jump}[1]{[\![#1]\!]}

\newcommand{\vect}[1]{\boldsymbol{#1}}

\title{A structure-preserving discontinuous Galerkin method for variable-density incompressible flow with a partitioned time integrator}
\author{Zonghao Bai \thanks{School of Mathematical Sciences, University of Science and Technology of China, Hefei, Anhui 230026, China. E-mail: baizonghao@mail.ustc.edu.cn.} \and Yinhua Xia\thanks{School of Mathematical Sciences, University of Science and Technology of China, Hefei, Anhui 230026, China. E-mail: yhxia@ustc.edu.cn. Research is supported in part by National Key R\&D Program of China No. 2022YFA1005202/2022YFA1005200, and NSFC grant No. 12671438.}  \and Yan Xu\thanks{Corresponding author. School of Mathematical Sciences, University of Science and Technology of China, Hefei, Anhui 230026, and Laoshan Laboratory, Qingdao 266237, P.R. China.  E-mail: yxu@ustc.edu.cn. The work of Yan Xu was partially supported by NSFC grant No. 12471347,  and Laoshan Laboratory (No. LSKJ202300305).
}}
\date{}

\begin{document}
\maketitle

\begin{abstract}
Preserving density bounds and mass conservation, together with a conservative
momentum discretization and exact incompressibility, is a central challenge in
high-order simulations of
variable-density incompressible flow.  We develop a structure-preserving
discontinuous Galerkin method that couples conservative density transport
with an \(H(\mathrm{div})\)-conforming mixed discretization of the momentum
 equation.  The method uses a partitioned time
 integrator: density is advanced explicitly by strong-stability-preserving
 Runge--Kutta stages with a
 bound-preserving limiter, while momentum is updated by a linearly implicit
 backward differentiation formula.
The transport velocities are extrapolated from previously computed mixed
velocities.  Because the discrete solenoidal space is linear, this
extrapolation preserves elementwise zero divergence and single-valued normal
traces, providing compatible transport velocities without an auxiliary
reconstruction or projection.  For source-free transport on triangular
meshes with impermeable or periodic boundaries, we prove bounds for density
cell averages and limiter control-point values under the stated spatial,
quadrature, and stagewise Courant--Friedrichs--Lewy conditions.  The
conservative fluxes and mean-preserving limiter yield global mass
conservation.  Numerical experiments on triangular and structured
rectangular meshes demonstrate the expected convergence rates, machine-level
divergence errors, and bound preservation, including in a low-density regime
with high density ratio.
\end{abstract}


\section{Introduction}

We consider the numerical approximation of variable-density incompressible
flow, which arises in buoyancy-driven transport and mixing
\citep{Li2021BP,LundgrenNazarov2024JCP}.  The governing equations couple
density transport to momentum evolution through an incompressible velocity
field.  For source-free transport with periodic or impermeable boundaries,
the density remains within its initial range and total mass is conserved.
Retaining these properties in a high-order approximation is important near
sharp density gradients, where numerical overshoots and undershoots can
affect both the inertial terms and gravitational forcing.  A suitable
discretization should therefore combine accurate density transport with
preservation of density bounds, mass conservation, and a velocity field
compatible with the incompressibility constraint.

Several approaches have been developed for this coupled problem.
Projection and splitting methods separate the velocity--pressure system
into subproblems, allowing a dedicated solver for each component
\citep{AlmgrenBellColellaHowellWelcome1998,BellMarcus1992,
GuermondQuartapelle2000,PyoShen2007}.
This decoupling simplifies the flow solve, while preservation of density
bounds remains dependent on the transport update.
Standard Taylor--Hood and MINI pairs provide
continuous velocities with single-valued normal traces, but their
incompressibility constraints do not generally yield elementwise
divergence-free velocities~\citep{John2017Divergence}.
Thus, satisfying the mixed incompressibility constraint alone need not
supply a transport velocity compatible with the bound-preserving (BP)
argument considered here.

For conservative discontinuous Galerkin (DG) density transport, the
Zhang--Shu limiter is bound-preserving (BP): it enforces bounds at selected
points by scaling the polynomial about its admissible cell average
\citep{ZhangShu2010,ZhangShu2011Survey}.
This correction preserves the cell average and hence conservation, but
requires the transport update to keep that average within the desired
bounds before limiting.
For variable-density flow, Li et al.~\citep{Li2021BP}
established a high-order
framework preserving both density bounds.  Their DG transport update uses
elementwise divergence-free velocities with single-valued normal traces,
obtained from a streamfunction construction alongside a continuous finite
element momentum solver; this construction entails an additional Poisson
solve at each transport stage.  The schemes of Li et al.~\citep{Li2024Energy} combine
energy stability with a lower density bound.  For transport, they use a
divergence-free projection into a Raviart--Thomas (RT) finite element
space~\citep{RaviartThomas1977}.  Their lower-bound guarantee
does not by itself ensure preservation of the initial upper density bound.
In these approaches, the local divergence-free property used in transport
is supplied by an auxiliary construction and need not hold for the
underlying momentum velocity.

Related semi-implicit DG treatments use a post-processed velocity in the
density equation and establish convergence for variable-density flow
\citep{LiQiuYang2022PostprocessedDG}.  Entropy-stable DG formulations provide
a complementary route to stability for variable-density incompressible
models, although the formulation in Manzanero et al.~\citep{ManzaneroRubioKopriva2020} uses
artificial compressibility rather than a strongly imposed incompressibility
constraint.  These works highlight the importance of making the transport
velocity and the conservation or stability mechanism compatible.

An alternative is to obtain the required velocity structure directly from
the momentum discretization.  \(H(\mathrm{div})\)-conforming spaces provide
single-valued normal traces; when the pressure space contains the velocity
divergences, the mixed constraint also yields elementwise zero divergence
\citep{BrezziDouglasMarini1985,BoffiBrezziFortin2013}.
Foundational divergence-free DG formulations for incompressible
Navier--Stokes equations were developed by Cockburn et al.~\citep{CockburnKanschatSchoetzau2007},
and related \(H(\mathrm{div})\)-conforming DG methods with semi-implicit
time-marching were analyzed by Han and Hou~\citep{HanHou2022HdivDG}.
Further \(H(\mathrm{div})\)-conforming and divergence-free DG formulations,
including unified comparisons with continuous Galerkin methods, are discussed
by Guzm\'an et al.~\citep{GuzmanShuSequeira2017}, Fu~\citep{Fu2019DivergenceFreeDG},
and Chen et al.~\citep{ChenLiDrapacaCimbala2020}.
The present paper also relates to the exactly divergence-free ultra-weak DG
developments in Xue and Xia~\citep{XueUWDG,XueYanXiaUWDG}.
These constructions provide the local velocity structure needed for
transport, but extending them to variable-density flow still requires a
density update that preserves the desired bounds.  Moreover, the transport
velocities must retain this structure at every time-integration stage.
This makes the coupling of density
transport, the mixed momentum solve, and velocity extrapolation a central
part of the method.

Conservation properties of variable-density discretizations have also been
studied through modified Galerkin formulations that preserve combinations
of mass, momentum, and energy identities~\citep{LundgrenNazarov2024JCP,
ZhangDongWang2023}.  The present method follows a different route: it keeps
the conservative momentum and density equations in flux form, while placing
the main structural emphasis on exact incompressibility and bound-preserving
density transport.

We address this coupling with a partitioned time integrator.  Its density
component is a strong-stability-preserving (SSP) Runge--Kutta method (SSPRK),
and its momentum component is a backward differentiation formula (BDF).  We
refer to this coupling as the SSPRK--BDF scheme.  On triangles, we use
Brezzi--Douglas--Marini (BDM) elements, while on rectangles we use
Raviart--Thomas (RT) elements~\citep{BrezziDouglasMarini1985,RaviartThomas1977}.
The corresponding pressure spaces are combined with an interior-penalty DG
viscosity discretization~\citep{Arnold1982IP}.
Density is advanced by SSPRK stages
\citep{GottliebSSP}, with a BP limiter applied at each stage, while momentum
is updated by BDF stages.
The key ingredient is the extrapolation of previously computed mixed
velocities: linearity of the discrete solenoidal space preserves their
elementwise zero divergence and single-valued normal traces.

The resulting coupling has two practical advantages.  First, the mixed
momentum solve supplies compatible velocities for the BP density update
without an additional streamfunction reconstruction or RT projection
during time stepping.  Second, the updated density and extrapolated
endpoint velocity are known in the momentum step, so the new velocity and
pressure are obtained from a linear Oseen-type mixed system.
The method thus combines explicit BP density transport with a linearly
implicit momentum update while retaining the required velocity structure.

Starting from divergence-free velocity histories, elementwise zero divergence is preserved at every solvable mixed step and under velocity extrapolation, independently of the prescribed density sequence. For source-free transport on triangular meshes with periodic or impermeable boundaries, it shows that the forward-Euler and SSPRK
bounds on cell averages under admissible starting histories, the stated
degree and quadrature assumptions, and a Courant--Friedrichs--Lewy (CFL)
condition for each actual extrapolated stage velocity.  The BP limiter
extends these bounds to its control points while preserving cell averages;
together with conservative fluxes, it ensures global mass conservation.
Accuracy tests on both mesh families assess observed convergence,
including in the low-density regime.
Classical benchmarks on triangular and rectangular meshes examine density
evolution, limiter effects, mass conservation, and divergence errors,
providing numerical evidence of the divergence-free and BP properties.
Total-energy histories serve as numerical diagnostics.


The rest of the paper is organized as follows. Section~\ref{sec:algfull} states the conservative model.
Section~\ref{sec:scheme} presents the structure-preserving discretization,
including the partitioned time integrator and its main properties.
Section~\ref{sec:numerical-results} reports verification and benchmark
calculations, followed by a discussion and conclusions.  The proof mechanism
is summarized in the propositions, with routine details connected to the
standard bound-preserving DG arguments cited there.

\section{Governing equations and conservation properties}
\label{sec:algfull}

The conservative formulation used below makes the conservation mechanism of
the discretization explicit.  In particular, the density equation is
discretized in flux form, while the momentum equation is written in terms of
the conservative momentum \(\rho\vect{u}\).  This formulation is equivalent
to the advective form for smooth incompressible solutions, but is better
suited to conservative DG transport and variable-density momentum updates.

Let \(\Omega\subset\mathbb R^2\) be a bounded connected polygonal domain
and \(T>0\).  For a constant dynamic viscosity coefficient \(\nu>0\),
the density \(\rho\), velocity
\(\vect{u}\), and pressure \(p\) satisfy the advective formulation
\begin{equation}
  \left\{
  \begin{alignedat}{3}
    \partial_t\rho+\vect{u}\cdot\nabla\rho
    &=0,
    &\qquad& \text{in }(0,T]\times\Omega,\\
    \rho(\partial_t\vect{u}+\vect{u}\cdot\nabla \vect{u})
    -\nu\Delta\vect{u}+\nabla p
    &=\vect{f},
    &\qquad& \text{in }(0,T]\times\Omega,\\
    \nabla\cdot \vect{u}
    &=0,
    &\qquad& \text{in }(0,T]\times\Omega .
  \end{alignedat}
  \right.
  \label{eq:algfull-nonconservative-system}
\end{equation}
For a smooth solution satisfying \(\nabla\cdot\vect{u}=0\), the advective
system \eqref{eq:algfull-nonconservative-system} is equivalent to
\begin{equation}
  \left\{
  \begin{alignedat}{3}
    \partial_t\rho+\nabla\cdot(\rho\vect{u})
    &=0,
    &\qquad& \text{in }(0,T]\times\Omega,\\
    \partial_t(\rho\vect{u})+\nabla\cdot(\rho\vect{u}\otimes \vect{u})
    -\nu\Delta\vect{u}+\nabla p
    &=\vect{f},
    &\qquad& \text{in }(0,T]\times\Omega,\\
    \nabla\cdot \vect{u}
    &=0,
    &\qquad& \text{in }(0,T]\times\Omega .
  \end{alignedat}
  \right.
  \label{eq:algfull-conservative-system}
\end{equation}
The equivalence follows from
the identities
\begin{equation}
  \partial_t\rho+\nabla\cdot(\rho\vect{u})
  =
  \partial_t\rho+\vect{u}\cdot\nabla\rho
  +\rho\nabla\cdot\vect{u},
  \label{eq:density-form-identity}
\end{equation}
and
\begin{align}
  \partial_t(\rho\vect{u})
  +\nabla\cdot(\rho\vect{u}\otimes\vect{u})
  &=
  \rho(\partial_t\vect{u}+\vect{u}\cdot\nabla\vect{u})
  +
  \vect{u}\bigl(\partial_t\rho+\nabla\cdot(\rho\vect{u})\bigr)
  \notag\\
  &=
  \rho(\partial_t\vect{u}+\vect{u}\cdot\nabla\vect{u}),
  \label{eq:algfull-product-identity}
\end{align}
where the last equality uses the density equation.
Let \(\vect{n}\) denote the outward unit normal and
\(\Gamma_D\subseteq\partial\Omega\) the boundary portion with prescribed
tangential velocity.  We write
\(\vect{u}_\tau=\vect{u}-(\vect{u}\cdot\vect{n})\vect{n}\).
The initial and boundary conditions are
\begin{equation}
  \left\{
  \begin{alignedat}{3}
    \rho(0,\vect{x})
    &=\rho_0(\vect{x}),
    &\qquad& \text{in }\Omega,\\
    \vect{u}(0,\vect{x})
    &=\vect{u}_0(\vect{x}),
    &\qquad& \text{in }\Omega,\\
    \vect{u}_\tau(t,\vect{x})
    &=\vect{u}_D(t,\vect{x}),
    &\qquad& \text{on }(0,T]\times\Gamma_D,
  \end{alignedat}
  \right.
  \label{eq:algfull-initial-boundary}
\end{equation}
where \(\vect{u}_D\cdot\vect{n}=0\) on \(\Gamma_D\).
All physical walls are impermeable:
\begin{equation}
  \vect{u}\cdot\vect{n}=0
  \qquad\text{on }\partial\Omega.
  \label{eq:impermeable-boundary-assumption}
\end{equation}
On \(\partial\Omega\setminus\Gamma_D\), the tangential viscous traction
vanishes in the natural boundary condition for the Laplacian.
In the periodic case, these wall
conditions are replaced by periodicity of the unknowns and data.

\section{Structure-preserving discretization}
\label{sec:scheme}

The method is designed to preserve four compatible structures: elementwise
zero divergence of the discrete velocity, single-valued normal velocity
traces, conservative density transport, and prescribed density bounds at the
limiter control points.  In this work, ``structure-preserving'' refers to
these properties; a fully discrete energy estimate is not claimed.

The discretization couples an explicit SSP/BP density update with a
linearly implicit mixed BDF momentum step.  Both subproblems use transport
velocities extrapolated from the discrete velocity history.  The resulting
partitioned time integrator is summarized before the spatial forms are
defined.

\subsection{Method overview}
\label{subsec:method-overview}

At time \(t_n\), the method proceeds as follows.  First, previously computed
mixed velocities are extrapolated to the endpoint and intermediate stages
required by the SSPRK method.  Second, the density is advanced through the three
SSPRK stages using conservative upwind DG fluxes; after each stage, the
Zhang--Shu limiter is applied and the limited state is used in the next
stage.  Finally, the updated density and endpoint-extrapolated velocity are
held fixed while one linearly implicit BDF3 mixed momentum solve produces
the new velocity and pressure.  Thus the density update is explicit and
bound-preserving, whereas the momentum update is implicit only through a
linear Oseen-type system.

The principal compatibility mechanism is velocity extrapolation within the
discrete solenoidal \(H(\mathrm{div})\) space.  Consequently, every stage
velocity has the elementwise zero-divergence and normal-continuity properties
needed by the density bound argument, without an auxiliary reconstruction or
projection.

\subsection{Notation and function spaces}

Let \(\mathcal T_h\) be a shape-regular triangulation of \(\Omega\),
and let \(\mathcal E_h^i\) and \(\mathcal E_h^\partial\) denote
the sets of interior and boundary faces, respectively.
We use \((\cdot,\cdot)_K\) and \(\langle\cdot,\cdot\rangle_F\)
for the usual \(L^2\) inner products on \(K\) and \(F\).

Let \(\vect{n}_K\) denote the outward unit normal to \(\partial K\).
On each interior face \(F=\partial K^-\cap\partial K^+\),
fix the labels \(K^-\) and \(K^+\), and choose \(\vect{n}_F\)
to point from \(K^-\) to \(K^+\).
Thus \(\vect{n}_{K^-}=\vect{n}_F\) and
\(\vect{n}_{K^+}=-\vect{n}_F\) on \(F\).
For a scalar, vector, or matrix field \(z\), denote its traces by
\(z^\pm=(z|_{K^\pm})|_F\), and define
\begin{equation}
  \avg{z}_F=\frac12(z^-+z^+),
  \qquad
  \jump{z}_F=z^--z^+.
  \label{eq:trace-notation}
\end{equation}
In the elementwise forms below, trial and test functions are
restricted to the current element \(K\), and their boundary
values are taken from the interior of \(K\).
Numerical fluxes are specified separately for each discrete operator.

We use the following function spaces:
\begin{align}
  H(\mathrm{div};\Omega)
  &:=\left\{
    \vect{v}\in[L^2(\Omega)]^2:
    \nabla\cdot\vect{v}\in L^2(\Omega)
  \right\},
  \label{eq:hdiv-space}\\
  H_0(\mathrm{div};\Omega)
  &:=\left\{
    \vect{v}\in H(\mathrm{div};\Omega):
    \left.\vect{v}\cdot\vect{n}\right|_{\partial\Omega}=0
  \right\},
  \label{eq:hdiv-zero-space}\\
  L_0^2(\Omega)
  &:=\left\{
    q\in L^2(\Omega):
    \int_\Omega q\,dx=0
  \right\}.
  \label{eq:l2-zero-space}
\end{align}

For polynomial degrees \(k_u\ge1\) and \(k_\rho\ge0\), the triangular
velocity, pressure, and density spaces are
\begin{align}
  V_h^{k_u}
  &=\{\vect{v}_h\in H_0(\mathrm{div};\Omega):
      \vect{v}_h|_K\in\mathrm{BDM}_{k_u}(K)=[\mathbb P_{k_u}(K)]^2
      \quad\forall K\in\mathcal T_h\},
  \label{eq:velocity-space}\\
  Q_h^{k_u-1}
  &=\{q_h\in L_0^2(\Omega):q_h|_K\in\mathbb P_{k_u-1}(K)\},
  \label{eq:pressure-space}\\
  W_h^{k_\rho}
  &=\{\varphi_h\in L^2(\Omega):
      \varphi_h|_K\in\mathbb P_{k_\rho}(K)\}.
  \label{eq:density-space}
\end{align}
For the triangular pair \(\mathrm{BDM}_{k_u}/\mathbb P_{k_u-1}\), the
following divergence inclusion and normal-continuity properties hold:
\begin{equation}
  \nabla\cdot V_h^{k_u}\subseteq Q_h^{k_u-1},
  \qquad
  \jump{\vect{v}_h\cdot\vect{n}_F}_F=0
  \quad(F\in\mathcal E_h^i).
  \label{eq:hdiv-structure}
\end{equation}
On a rectangle \(K\), let \(\mathbb Q_{a,b}(K)\) denote polynomials
of degree at most \(a\) and \(b\) in the two coordinate directions,
and set \(\mathbb Q_k(K)=\mathbb Q_{k,k}(K)\).  We use
\[
  \mathrm{RT}_{k_u}(K)
  =\mathbb Q_{k_u+1,k_u}(K)\times\mathbb Q_{k_u,k_u+1}(K)
\]
for velocity, discontinuous \(\mathbb Q_{k_u}(K)\) for pressure
with global zero mean, and discontinuous \(\mathbb Q_{k_\rho}(K)\)
for density.  The \(\mathrm{RT}_{k_u}/\mathbb Q_{k_u}\) pair has the
same divergence inclusion and normal-continuity properties
\citep{RaviartThomas1977,BoffiBrezziFortin2013}.  The analysis below is written
for triangular BDM elements.
For periodic boundary conditions, the zero normal boundary condition in
\eqref{eq:velocity-space} is replaced by matching normal traces on paired
periodic faces, with each pair treated as a single interior interface in
\(\mathcal E_h^i\) for the traces and discrete forms.

Define the divergence-free velocity subspace by
\begin{equation}
  Z_h:=\left\{\vect{v}_h\in V_h^{k_u}:
    \left.(\nabla\cdot\vect{v}_h)\right|_K=0
    \quad\forall K\in\mathcal T_h\right\}.
  \label{eq:discrete-solenoidal-space}
\end{equation}
Both the computed momentum velocity \(\vect{u}_h\) and the prescribed
transport velocity \(\vect{w}_h\) are represented in \(V_h^{k_u}\).
The notation \(Z_h\) records the additional elementwise divergence-free
property used in the analysis.

\subsection{Conservative DG density transport}
\label{subsec:density-equation}

For a prescribed transport velocity \(\vect{w}\), the density subproblem is
\begin{equation}
  \partial_t\rho+\nabla\cdot(\rho\vect{w})=0,
  \qquad
  \nabla\cdot\vect{w}=0.
  \label{eq:density-subproblem}
\end{equation}
For a prescribed discrete transport velocity \(\vect{w}_h\in V_h^{k_u}\)
satisfying \(\nabla\cdot\vect{w}_h=0\) on each element, the density
approximation is sought in \(W_h^{k_\rho}\).
Proposition~\ref{prop:divergence-independent-density} verifies this velocity
condition for the momentum update and its extrapolations.

\subsubsection{Spatial discretization}
\label{subsubsec:density-space}

Let \(\vect{w}_h\in V_h^{k_u}\) be the prescribed transport velocity above.
Its normal continuity across interior faces allows us to define
\begin{equation}
  a_F=\vect{w}_h^-\cdot\vect{n}_F
     =\vect{w}_h^+\cdot\vect{n}_F,
  \qquad
  a_{K,F}=(\vect{w}_h|_K)\cdot\vect{n}_K,
  \label{eq:normal-speed}
\end{equation}
where \(a_F\) is defined on interior faces and \(a_{K,F}\) on each
face of \(K\).  For a scalar or vector quantity \(z\), define its
upwind trace pointwise by
\begin{equation}
  z_F^{\mathrm{up}}
  :=\begin{cases}
    z^-,&a_F\ge0,\\
    z^+,&a_F<0,
  \end{cases}
  \qquad F\in\mathcal E_h^i.
  \label{eq:transport-upwind-trace}
\end{equation}
The outward numerical density flux on \(F\subset\partial K\) is
\begin{equation}
  \widehat f_{\rho,K,F}
  :=\begin{cases}
    a_{K,F}\rho_{h,F}^{\mathrm{up}},
      &F\in\mathcal E_h^i,\\
    0,
      &F\in\mathcal E_h^\partial.
  \end{cases}
  \label{eq:density-upwind-flux}
\end{equation}
Paired periodic faces are treated as interior interfaces.
The physical boundary flux vanishes by impermeability.
On an interior interface, \(a_{K^-,F}=a_F=-a_{K^+,F}\), so the
two elements have equal and opposite numerical fluxes.

For \(\rho_h,\varphi_h\in W_h^{k_\rho}\), define the density residual
\begin{equation}
  \mathcal R_h(\rho_h,\vect{w}_h;\varphi_h)
  :=
  \sum_{K\in\mathcal T_h}
  \left[
    (\rho_h\vect{w}_h,\nabla\varphi_h)_K
    -\sum_{F\subset\partial K}
     \langle\widehat f_{\rho,K,F},
     \varphi_h\rangle_{F,h}
  \right].
  \label{eq:density-residual}
\end{equation}
Here \(\langle\cdot,\cdot\rangle_{F,h}\) denotes the face quadrature
pairing defined in \eqref{eq:density-face-quadrature}.
The semi-discrete problem is to find
\(\rho_h(t)\in W_h^{k_\rho}\) such that
\begin{equation}
  (\partial_t\rho_h,\varphi_h)_\Omega
  =\mathcal R_h(\rho_h,\vect{w}_h;\varphi_h)
  \qquad\forall\varphi_h\in W_h^{k_\rho}.
  \label{eq:density-semidiscrete}
\end{equation}
Taking \(\varphi_h=\mathbf1_K\), the indicator function of \(K\), gives
\begin{equation}
  |K|\frac{d\bar\rho_K}{dt}
  =-\sum_{F\subset\partial K}
       \langle\widehat f_{\rho,K,F},1\rangle_{F,h},
  \qquad
  \bar\rho_K=\frac1{|K|}\int_K\rho_h\,dx.
  \label{eq:density-average-balance}
\end{equation}
Summation over \(K\) yields the conservative DG mass balance by cancellation
of the interior fluxes \citep{CockburnShu2001}.

\subsubsection{Density time discretization}
\label{subsubsec:density-time}

Let \(t_n=n\Delta t\), with a uniform time step \(\Delta t>0\).
The first-order SSP method is the forward-Euler update
\begin{equation}
  (\rho_h^{n+1},\varphi_h)_\Omega
  =
  (\rho_h^n,\varphi_h)_\Omega
  +\Delta t\,\mathcal R_h(\rho_h^n,\vect{w}_h^n;\varphi_h)
  \qquad\forall\varphi_h\in W_h^{k_\rho}.
  \label{eq:density-ssp1}
\end{equation}
The second-order SSP method is
\begin{align}
  (\rho_h^{(1)},\varphi_h)_\Omega
  &=
  (\rho_h^n,\varphi_h)_\Omega
  +\Delta t\,\mathcal R_h(\rho_h^n,\vect{w}_h^n;\varphi_h),
  \notag\\
  (\rho_h^{n+1},\varphi_h)_\Omega
  &=
  \frac12(\rho_h^n,\varphi_h)_\Omega
  +\frac12\left[
    (\rho_h^{(1)},\varphi_h)_\Omega
    +\Delta t\,\mathcal R_h(\rho_h^{(1)},\vect{w}_h^{n+1};\varphi_h)
  \right]
  \qquad\forall\varphi_h\in W_h^{k_\rho}.
  \label{eq:density-ssp2}
\end{align}
The third-order SSP Runge--Kutta method is
\begin{align}
  (\rho_h^{(1)},\varphi_h)_\Omega
  &=
  (\rho_h^n,\varphi_h)_\Omega
  +\Delta t\,\mathcal R_h(\rho_h^n,\vect{w}_h^n;\varphi_h),
  \notag\\
  (\rho_h^{(2)},\varphi_h)_\Omega
  &=
  \frac34(\rho_h^n,\varphi_h)_\Omega
  +\frac14\left[
    (\rho_h^{(1)},\varphi_h)_\Omega
    +\Delta t\,\mathcal R_h(\rho_h^{(1)},\vect{w}_h^{n+1};\varphi_h)
  \right],
  \notag\\
  (\rho_h^{n+1},\varphi_h)_\Omega
  &=
  \frac13(\rho_h^n,\varphi_h)_\Omega
  +\frac23\left[
    (\rho_h^{(2)},\varphi_h)_\Omega
    +\Delta t\,\mathcal R_h(\rho_h^{(2)},\vect{w}_h^{n+1/2};\varphi_h)
  \right].
  \label{eq:density-ssprk3}
\end{align}
These identities hold for every \(\varphi_h\in W_h^{k_\rho}\).
The stage velocities satisfy \(\vect{w}_h^n=\vect{u}_h^n\).  Their endpoint
extrapolations, together with the third-order midpoint extrapolation, are
given by
\begin{align}
  r=1:\quad&
  \vect{w}_h^{n+1}=\vect{u}_h^n,
  \notag\\
  r=2:\quad&
  \vect{w}_h^{n+1}=2\vect{u}_h^n-\vect{u}_h^{n-1},
  \notag\\
  r=3:\quad&
  \vect{w}_h^{n+1}
  =3\vect{u}_h^n-3\vect{u}_h^{n-1}+\vect{u}_h^{n-2}.
  \notag\\
  &
  \vect{w}_h^{n+1/2}
  =\frac{15}{8}\vect{u}_h^n-\frac54\vect{u}_h^{n-1}
   +\frac38\vect{u}_h^{n-2}.
  \label{eq:density-endpoint-extrapolation}
\end{align}
Equations~\eqref{eq:density-ssp1}--\eqref{eq:density-ssprk3} define the
SSP updates before limiting.  In the BP scheme, each completed stage is
replaced by its BP projection; all subsequent residuals and convex
combinations use the limited state.

\subsubsection{Bound-preserving density update}
\label{subsubsec:density-bp}

Consider \eqref{eq:density-subproblem} with the impermeable wall condition
\eqref{eq:impermeable-boundary-assumption} or paired periodic boundaries.
Let \(x_{F,\beta}\), \(\beta=1,\ldots,k_\rho+1\), be the mapped Gauss
points on a face \(F\), with weights \(\varpi_\beta>0\) satisfying
\(\sum_{\beta=1}^{k_\rho+1}\varpi_\beta=1\).
The face quadrature pairing used in the density discretization is
\begin{equation}
  \langle g,\psi\rangle_{F,h}
  :=|F|\sum_{\beta=1}^{k_\rho+1}\varpi_\beta
       g(x_{F,\beta})\psi(x_{F,\beta}).
  \label{eq:density-face-quadrature}
\end{equation}
For constant-state consistency of the density discretization, impose
\begin{equation}
  k_u\le k_\rho+1.
  \label{eq:density-degree-condition}
\end{equation}
The face Gauss rule then integrates
\((\vect{w}_h\cdot\vect{n}_K)\varphi_h\) exactly for
\(\varphi_h\in W_h^{k_\rho}\), so the density residual vanishes
for constant density and \(\vect{w}_h\in Z_h\).
Let \(S_K\) consist of the nodes in the positive triangular quadrature of
Zhang, Xia, and Shu~\citep[Section~3]{ZhangXiaShu2012JSC}, including the face nodes of
\eqref{eq:density-face-quadrature}.
For fixed density bounds \(0<m<M\), define
\begin{equation}
  \mathcal A_h(m,M)
  :=\left\{q_h\in W_h^{k_\rho}:
    m\le q_h|_K(\xi)\le M
    \quad\forall K\in\mathcal T_h,\ \forall\xi\in S_K\right\}.
  \label{eq:discrete-admissible-set}
\end{equation}
Positive quadrature gives \(\bar q_K\in[m,M]\) for
\(q_h\in\mathcal A_h(m,M)\).  Here bound preservation refers to
cell averages and values at the limiter control points \(S_K\).

\paragraph{Bound-preserving limiter.}
For \(q_h\in W_h^{k_\rho}\) with \(\bar q_K\in[m,M]\), set
\(q_K=q_h|_K\) and
\begin{equation}
  \bar q_K=\frac1{|K|}\int_Kq_h\,dx,
  \qquad
  q_K^{\min}=\min_{\xi\in S_K}q_K(\xi),
  \qquad
  q_K^{\max}=\max_{\xi\in S_K}q_K(\xi),
  \label{eq:bp-stage-extrema}
\end{equation}
and
\begin{equation}
  \theta_K
  =\min\left\{
    1,
    \frac{M-\bar q_K}{q_K^{\max}-\bar q_K},
    \frac{\bar q_K-m}{\bar q_K-q_K^{\min}}
  \right\},
  \label{eq:bp-theta}
\end{equation}
where a quotient with zero denominator is assigned the value one.  The
Zhang--Shu scaling \citep{ZhangShu2010,ZhangShu2011Survey} is
\begin{equation}
  (\Pi_h^{\mathrm{BP}}q_h)|_K
  =\bar q_K+\theta_K(q_K-\bar q_K).
  \label{eq:bp-scaling}
\end{equation}

The standard scaling properties
\citep{ZhangShu2010,ZhangShu2011Survey} give \(0\le\theta_K\le1\) and
\begin{equation}
  \Pi_h^{\mathrm{BP}}q_h\in\mathcal A_h(m,M),
  \qquad
  \overline{\Pi_h^{\mathrm{BP}}q_h}_K=\bar q_K.
  \label{eq:bp-projection-properties}
\end{equation}
The SSP updates \eqref{eq:density-ssp1}--\eqref{eq:density-ssprk3}
and the mean-preserving BP scaling preserve the global mass balance
obtained by summing \eqref{eq:density-average-balance} over all elements.

Let \(F_i\), \(i=1,2,3\), be the edges of \(K\), with outward normals
\(\vect{n}_{K,i}\), and set
\begin{equation}
  \alpha_K(\vect{w}_h)
  =\max_{\substack{1\le i\le3\\1\le\beta\le k_\rho+1}}
    \left|(\vect{w}_h|_K)(x_{i,\beta})\cdot\vect{n}_{K,i}\right|.
  \label{eq:bp-common-speed}
\end{equation}
Let \(\widehat\omega_1\) denote the endpoint weight of the normalized
\(N_q\)-point Gauss--Lobatto rule, where \(N_q\ge2\) is the smallest
integer satisfying \(2N_q-3\ge k_\rho\).  This is the quadrature weight
in the one-dimensional edge parameter of the positive triangular
quadrature.  Thus, although the cell quadrature is two-dimensional, the
corresponding weight at the \(\beta\)-th Gauss point on edge \(F_i\) is
\(\omega_{i,\beta}=\frac23\widehat\omega_1\varpi_\beta\), which is the
individual edge-node weight. 
Condition~\eqref{eq:density-degree-condition} implies
\(k_u\le2k_\rho+1\), as required for the following cell-average bounds.

\begin{proposition}[Forward-Euler cell-average bound]
\label{prop:forward-euler-average-bound}
Let \(\rho_h^n\in\mathcal A_h(m,M)\) and \(\vect{w}_h\in Z_h\).
Denote by \(\rho_h^{\mathrm{FE}}\) the forward-Euler update
\eqref{eq:density-ssp1} with this velocity.  If
\begin{equation}
  \frac{\Delta t\,\alpha_K(\vect{w}_h)}{|K|}
  \sum_{i=1}^3|F_i|
  \le\frac23\widehat\omega_1
  \qquad\forall K\in\mathcal T_h,
  \label{eq:bp-cfl}
\end{equation}
then
\begin{equation}
  m\le\bar\rho_K^{\mathrm{FE}}\le M
  \qquad\forall K\in\mathcal T_h.
  \label{eq:forward-euler-average-bound}
\end{equation}
\end{proposition}

The result follows from the standard positive-quadrature convex-combination
argument for bound-preserving DG transport.  Under the stated CFL condition,
the coefficients in that decomposition are nonnegative.  Since the transport
velocity is elementwise divergence-free, their sum is one.  The cell and
upwind face values therefore lie in the admissible interval, which gives the
claimed bound.  See Zhang and Shu~\citep{ZhangShu2010,ZhangShu2011Survey}
and the triangular construction of Zhang, Xia, and Shu~\citep{ZhangXiaShu2012JSC}
for the corresponding positive-quadrature decompositions.


\begin{proposition}[SSPRK3 bound preservation]
\label{prop:ssprk3-bp}
Let \(\rho_h^n\in\mathcal A_h(m,M)\) and
\(\vect{u}_h^{n-j}\in Z_h\), \(j=0,1,2\).
For the stage velocities defined in
\eqref{eq:density-endpoint-extrapolation},
assume
\begin{equation}
  \frac{\Delta t\,\alpha_K(\vect{w}_h^{n+c})}{|K|}
  \sum_{i=1}^3|F_i|
  \le\frac23\widehat\omega_1
  \qquad\forall K\in\mathcal T_h,
  \quad c\in\left\{0,\frac12,1\right\}.
  \label{eq:ssprk3-stage-cfl}
\end{equation}
Then, after applying \eqref{eq:bp-scaling} at each stage, the provisional
cell averages lie in \([m,M]\).  The limited states satisfy
\begin{equation}
  \rho_h^{(1)},\ \rho_h^{(2)},\ \rho_h^{n+1}
  \in\mathcal A_h(m,M),
  \qquad
  \bar\rho_K^{(1)},\ \bar\rho_K^{(2)},\ \bar\rho_K^{n+1}
  \in[m,M].
  \label{eq:ssprk3-average-bound-chain}
\end{equation}
\end{proposition}

For the first-order scheme, set
\(\rho_h^{n+1}=\Pi_h^{\mathrm{BP}}\rho_h^{\mathrm{FE}}\).  The third-order
SSPRK method is a convex combination of forward-Euler updates.  Moreover, the
extrapolated stage velocities remain in the linear solenoidal space, so
Proposition~\ref{prop:forward-euler-average-bound} applies at each stage under
the stagewise CFL condition.  The mean-preserving Zhang--Shu limiter then
enforces the control-point bounds without changing cell averages; see Zhang
and Shu~\citep{ZhangShu2011Survey}.


The same argument applies to SSP2 in \eqref{eq:density-ssp2}, provided
\(\vect{w}_h^n,\vect{w}_h^{n+1}\in Z_h\) and \eqref{eq:bp-cfl}
holds for both stage velocities.
Condition~\eqref{eq:ssprk3-stage-cfl} uses the actual extrapolated stage
velocities, whose normal speeds can exceed those of the velocity history.
Although the density stages retain the classical SSPRK3 convex-combination
structure \citep{GottliebSSP}, the present BP argument must also account for
velocity extrapolation.  The extrapolation coefficients can be negative,
so the proof separately uses the preservation of \(Z_h\) under extrapolation
and the CFL condition for each extrapolated stage velocity.  These
observations are the standard SSP extension of the forward-Euler bound
argument; see Zhang and Shu~\citep{ZhangShu2011Survey}.

\begin{remark}[Tensor-product quadrilateral extension]

For the RT/\(\mathbb Q\) pair on rectangular meshes, use the
scalar maximum-principle control set of
\citep[Section~3.2]{ZhangShu2010}, including the Gauss/Gauss and both
mixed Gauss--Lobatto/Gauss products.  The positive cell-average
decomposition on the mixed products is the standard tensor-product
construction; see
\citep{ZhangShu2010,ZhangShu2010EulerRect,ZhangShu2011Survey}.
With four edges in place of three and \(\widehat\omega_1/2\) in place of
\(2\widehat\omega_1/3\) in the stage CFL conditions, the same argument
proves bounds on cell averages and limiter control-point values.

The structure-preserving arguments are not tied to a particular
\(H(\operatorname{div})\)-conforming velocity element.  They apply equally
to BDM, Raviart--Thomas (RT), and their tensor-product counterparts on
Cartesian meshes, provided that the selected pair has single-valued normal
traces, supports the required discrete divergence constraint, and has the
polynomial approximation properties used by the method.  Thus, the choice
between BDM and RT affects the approximation space, degrees of freedom, and
implementation details, but not the mass-conservation, bound-preservation,
or divergence-preservation mechanisms established here.
\end{remark}

\subsection{Mixed momentum discretization}
\label{subsec:momentum-equation}

For prescribed density \(\rho\) and Oseen velocity \(\vect{w}\), the
linearized conservative momentum--pressure subsystem is
\begin{equation}
  \partial_t(\rho\vect{u})
  +\nabla\cdot(\rho\vect{u}\otimes\vect{w})
  -\nu\Delta\vect{u}+\nabla p
  =\vect{f},
  \qquad
  \nabla\cdot\vect{u}=0.
  \label{eq:momentum-subproblem}
\end{equation}
The discrete unknowns are
\((\vect{u}_h(t),p_h(t))\in V_h^{k_u}\times Q_h^{k_u-1}\), with
\(\rho_h\) and \(\vect{w}_h\) prescribed.

\subsubsection{Spatial discretization}
\label{subsubsec:momentum-space}

For \(\rho_h\in W_h^{k_\rho}\) and \(\vect{w}_h\in V_h^{k_u}\) with
\(\nabla\cdot\vect{w}_h=0\) on each element, and for
\(\vect{u}_h,\vect{v}_h\in V_h^{k_u}\), set
\(\vect{m}_h=\rho_h\vect{u}_h\).  We use
\((\vect{a}\otimes\vect{b})_{ij}=a_i b_j\) and
\((\nabla\vect{v})_{ij}=\partial_j v_i\).

\textbf{Convection.}
With \(a_{K,F}\) and the upwind trace defined in
\eqref{eq:normal-speed}--\eqref{eq:transport-upwind-trace}, let
\begin{equation}
  \widehat{\vect F}_{m,K,F}
  :=\begin{cases}
    a_{K,F}\vect{m}_{h,F}^{\mathrm{up}},&F\in\mathcal E_h^i,\\
    \vect{0},&F\in\mathcal E_h^\partial,
  \end{cases}
  \qquad F\subset\partial K.
  \label{eq:momentum-upwind-flux}
\end{equation}
The conservative convection form is
\begin{equation}
  \mathcal C_h(\rho_h,\vect{u}_h,\vect{w}_h;\vect{v}_h)
  :=\sum_{K\in\mathcal T_h}
  \left[
    -(\rho_h\vect{u}_h\otimes\vect{w}_h,\nabla\vect{v}_h)_K
    +\sum_{F\subset\partial K}
     \langle\widehat{\vect F}_{m,K,F},
     \vect{v}_h\rangle_F
  \right].
  \label{eq:momentum-convection}
\end{equation}
On each interior face,
\begin{equation}
  \widehat{\vect F}_{m,K^-,F}
  =a_F\vect{m}_{h,F}^{\mathrm{up}}
  =-\widehat{\vect F}_{m,K^+,F},
  \qquad F\in\mathcal E_h^i.
  \label{eq:momentum-flux-conservation}
\end{equation}

\textbf{Viscosity.}
For the symmetric interior-penalty discretization \citep{Arnold1982IP},
define
\begin{equation}
  \mathcal E_h^D
  :=\{F\in\mathcal E_h^\partial:F\subset\Gamma_D\}.
  \label{eq:dirichlet-face-set}
\end{equation}
Let \(h_F\) denote the local mesh size associated with \(F\).
The penalty coefficient is
\begin{equation}
  \tau_F:=\frac{C_{\mathrm{IP}}k_u^2}{h_F},
  \qquad C_{\mathrm{IP}}>0.
  \label{eq:ipdg-penalty}
\end{equation}
The constant \(C_{\mathrm{IP}}\) is chosen sufficiently large for
coercivity; all reported computations use \(C_{\mathrm{IP}}=6\).

For the bilinear form, define the numerical traces
\begin{equation}
  \widehat{\vect{u}}_h
  :=\begin{cases}
    \avg{\vect{u}_h}_F,&F\in\mathcal E_h^i,\\
    \vect{0},&F\in\mathcal E_h^D,
  \end{cases}
  \label{eq:ipdg-velocity-flux}
\end{equation}
\begin{equation}
  \widehat{\nabla\vect{u}_h}
  :=\begin{cases}
    \avg{\nabla\vect{u}_h}_F
      -\tau_F\jump{\vect{u}_h}_F\otimes\vect{n}_F,
      &F\in\mathcal E_h^i,\\
    \nabla\vect{u}_h-\tau_F\vect{u}_h\otimes\vect{n}_K,
      &F\in\mathcal E_h^D.
  \end{cases}
  \label{eq:ipdg-gradient-flux}
\end{equation}
The viscous bilinear form is
\begin{align}
  \mathcal A_h^{\mathrm{IP}}(\vect{u}_h;\vect{v}_h)
  &:=
  \sum_{K\in\mathcal T_h}\Biggl[
    (\nabla\vect{u}_h,\nabla\vect{v}_h)_K
  \notag\\
  &\qquad
    -\sum_{\substack{F\subset\partial K\\
                     F\in\mathcal E_h^i\cup\mathcal E_h^D}}
      \Bigl[
        \left\langle\widehat{\nabla\vect{u}_h}\vect{n}_K,
                     \vect{v}_h\right\rangle_F
        +\left\langle\nabla\vect{v}_h\vect{n}_K,
                     \vect{u}_h-\widehat{\vect{u}}_h\right\rangle_F
      \Bigr]
    \Biggr].
  \label{eq:ipdg-form}
\end{align}
Tangential Dirichlet data enter through the Nitsche functional
\citep{Nitsche1971},
\begin{equation}
  \mathcal L_{D,h}^{\mathrm{IP}}(\vect{u}_D;\vect{v}_h)
  :=\sum_{K\in\mathcal T_h}
    \sum_{\substack{F\subset\partial K\\F\in\mathcal E_h^D}}
    \left[
      -\langle\nabla\vect{v}_h\vect{n}_K,\vect{u}_D\rangle_F
      +\langle\tau_F\vect{u}_D,\vect{v}_h\rangle_F
    \right].
  \label{eq:ipdg-dirichlet-functional}
\end{equation}
Faces in \(\mathcal E_h^\partial\setminus\mathcal E_h^D\) carry the
natural tangential traction condition and contribute no boundary term.
On interior faces, set
\(\vect{u}_{h,\tau}=\vect{u}_h-(\vect{u}_h\cdot\vect{n}_F)\vect{n}_F\).
Normal conformity and the wall conditions give
\begin{equation}
  \begin{aligned}
    \jump{\vect{u}_h}_F
      &=\jump{\vect{u}_{h,\tau}}_F,
      &&F\in\mathcal E_h^i,\\
    (\vect{u}_h-\vect{u}_D)\cdot\vect{n}_K
      &=\vect{v}_h\cdot\vect{n}_K=0,
      &&F\in\mathcal E_h^D.
  \end{aligned}
  \label{eq:ipdg-tangential-jump}
\end{equation}
Thus the interior and boundary penalties act only on tangential traces.

\begin{remark}[Transient penalty]
\label{rem:algfull-transient-boundary}
Let
\(\rho_{h,F}=(\rho_h|_K)|_F\) for
\(F\subset\partial K\), \(F\in\mathcal E_h^D\).  The optional
transient penalty in the fully discrete momentum equation is
\begin{equation}
  \mathcal B_{\partial,h}^{\rho}
  (\rho_h,\vect{z};\vect{v}_h)
  :=\sum_{K\in\mathcal T_h}
    \sum_{\substack{F\subset\partial K\\F\in\mathcal E_h^D}}
    \left\langle
      \eta_{\rho,F}\vect{z}_\tau,
      \vect{v}_{h,\tau}
    \right\rangle_F.
  \label{eq:algfull-B-rho}
\end{equation}
Using the transient Nitsche scaling of
Schott et al.~\citep[Eq.~(35)]{Schott2015Nitsche}, set
\begin{equation}
  \eta_{\rho,F}
  =\gamma_D\frac{\rho_{h,F}h_F}{\Delta t},
  \qquad \gamma_D\ge0.
  \label{eq:algfull-eta-rho}
\end{equation}
The density trace is evaluated at \(t_{n+1}\).  For prescribed tangential
data, the penalty acts on \(\vect{z}=\vect{u}_h^{n+1}-\vect{u}_D^{n+1}\).
\end{remark}

\textbf{Pressure--velocity coupling.}
Define
\begin{equation}
	\mathcal G_h(q_h;\vect{v}_h)
	:=\sum_{K\in\mathcal T_h}(q_h,\nabla\cdot\vect{v}_h)_K,
	\qquad
	(q_h,\vect{v}_h)\in Q_h^{k_u-1}\times V_h^{k_u}.
	\label{eq:pressure-divergence-forms}
\end{equation}

With \(\rho_h\) and \(\vect{w}_h\) prescribed, the semi-discrete mixed problem seeks
\((\vect{u}_h(t),p_h(t))\in V_h^{k_u}\times Q_h^{k_u-1}\) such that
\begin{align}
  (\partial_t(\rho_h\vect{u}_h),\vect{v}_h)_\Omega
  &+\mathcal C_h(\rho_h,\vect{u}_h,\vect{w}_h;\vect{v}_h)
   +\nu\mathcal A_h^{\mathrm{IP}}(\vect{u}_h;\vect{v}_h)
    -\mathcal G_h(p_h;\vect{v}_h)
  \notag\\
  &\qquad
   =(\vect{f},\vect{v}_h)_\Omega
   +\nu\mathcal L_{D,h}^{\mathrm{IP}}(\vect{u}_D;\vect{v}_h)
  \quad\forall\vect{v}_h\in V_h^{k_u},
  \label{eq:momentum-semidiscrete}\\
  \mathcal G_h(q_h;\vect{u}_h)&=0
  \quad\forall q_h\in Q_h^{k_u-1}.
  \label{eq:divergence-semidiscrete}
\end{align}
The time-step-dependent penalty \(\mathcal B_{\partial,h}^{\rho}\) enters
only the fully discrete formulation.

\subsubsection{Momentum time discretization}
\label{subsubsec:momentum-time}

Let \(\vect{m}_h^j=\rho_h^j\vect{u}_h^j\), with \(\rho_h^{n+1}\)
obtained from the SSP/BP update.  The Oseen velocity
\(\vect{w}_h^{n+1}\) is the order-\(r\) endpoint extrapolation in
\eqref{eq:density-endpoint-extrapolation}.  For \(r\in\{1,2,3\}\), define
\begin{equation}
  \delta_r\vect{m}_h^{n+1}
  :=\frac1{\Delta t}\sum_{j=0}^r d_j^{(r)}\vect{m}_h^{n+1-j},
  \label{eq:bdf-difference}
\end{equation}
where
\begin{equation}
  \begin{aligned}
    (d_0^{(1)},d_1^{(1)})&=(1,-1),\\
    (d_0^{(2)},d_1^{(2)},d_2^{(2)})&=\left(\frac32,-2,\frac12\right),\\
    (d_0^{(3)},d_1^{(3)},d_2^{(3)},d_3^{(3)})
      &=\left(\frac{11}{6},-3,\frac32,-\frac13\right).
  \end{aligned}
  \label{eq:bdf-coefficients}
\end{equation}
The BDF mixed step seeks
\((\vect{u}_h^{n+1},p_h^{n+1})\in V_h^{k_u}\times Q_h^{k_u-1}\)
such that
\begin{align}
  &(\delta_r\vect{m}_h^{n+1},\vect{v}_h)_\Omega
   +\mathcal C_h(\rho_h^{n+1},\vect{u}_h^{n+1},
                  \vect{w}_h^{n+1};\vect{v}_h)
   +\nu\mathcal A_h^{\mathrm{IP}}(\vect{u}_h^{n+1};\vect{v}_h)
  \notag\\
  &\qquad
   +\mathcal B_{\partial,h}^{\rho}
      (\rho_h^{n+1},\vect{u}_h^{n+1}-\vect{u}_D^{n+1};\vect{v}_h)
   -\mathcal G_h(p_h^{n+1};\vect{v}_h)
   =(\vect{f}^{n+1},\vect{v}_h)_\Omega
    +\nu\mathcal L_{D,h}^{\mathrm{IP}}
       (\vect{u}_D^{n+1};\vect{v}_h)
  \label{eq:bdf-momentum}
\end{align}
for every \(\vect{v}_h\in V_h^{k_u}\), together with
\begin{equation}
  \mathcal G_h(q_h;\vect{u}_h^{n+1})=0
  \qquad\forall q_h\in Q_h^{k_u-1}.
  \label{eq:bdf-divergence}
\end{equation}
The coefficients \(\rho_h^{n+1}\) and \(\vect{w}_h^{n+1}\) are fixed;
hence the system is linear in the new velocity and pressure.
For an analysis of extrapolation-based linearly implicit BDF methods for
quasilinear parabolic equations, see Akrivis and Lubich~\citep{AkrivisLubich2015}.

Moreover, since the discrete solenoidal space \(Z_h\) is linear, any
extrapolation of previously computed velocities in \(Z_h\) remains
elementwise divergence-free and retains single-valued normal traces.  The
following proposition records this property for the BDF update and all
transport-velocity extrapolations used by the method.

\begin{proposition}[Divergence-free preservation under time discretization]
\label{prop:divergence-independent-density}
Let \(r\in\{1,2,3\}\), \(N\ge r\), and
\(\vect{u}_h^j\in Z_h\), \(j=0,\ldots,r-1\).
For any prescribed density sequence
\(\{\rho_h^j\}_{j=0}^N\subset W_h^{k_\rho}\), suppose that
\eqref{eq:bdf-momentum}--\eqref{eq:bdf-divergence} admits an exact
solution \((\vect{u}_h^{n+1},p_h^{n+1})\in V_h^{k_u}\times Q_h^{k_u-1}\)
at each step \(n=r-1,\ldots,N-1\).  Then, in \eqref{eq:density-endpoint-extrapolation} for every
\(K\in\mathcal T_h\),
\begin{equation}
  \begin{aligned}
    \nabla\cdot\vect{w}_h^{n+c}&=0,
      &&n=r-1,\ldots,N-1,\quad c\in\mathcal C_r,
  \end{aligned}
  \label{eq:divergence-after-step}
\end{equation}
where \(\mathcal C_1=\mathcal C_2=\{0,1\}\) and
\(\mathcal C_3=\{0,\tfrac12,1\}\).
\end{proposition}



For source-free transport, the positive quadrature and the stagewise CFL
conditions make each forward-Euler cell average a convex combination of
admissible cell and upwind face values.  The SSPRK3 stages inherit this
property through their convex-combination representation.  The
mean-preserving BP scaling then places the density at the prescribed bounds
on the control set without changing any cell average.  The rectangular RT construction is
the corresponding tensor-product extension used in the numerical tests.

\subsection{SSPRK--BDF main loop}
\label{subsec:complete-algorithm}

To initialize the third-order main loop, the states at \(t_1\) and \(t_2\)
can be generated by a one-step method of at least third-order accuracy
applied to the coupled semi-discrete system.  For source-free transport,
the startup must preserve density admissibility and produce velocities
satisfying the mixed incompressibility constraint, so that the resulting
histories satisfy the assumptions below.

For \(n\ge2\), suppose that \(\rho_h^n\in\mathcal A_h(m,M)\) and
\(\vect{u}_h^{n-j}\in Z_h\), \(j=0,1,2\).
Given the history
\[
  (\rho_h^{n-j},\vect{u}_h^{n-j},p_h^{n-j}),
  \qquad j=0,1,2,
\]
and a uniform time step satisfying \eqref{eq:ssprk3-stage-cfl}, the
SSPRK3--BDF3 update is defined by the following three steps.
\begin{enumerate}
  \item Set \(\vect{w}_h^n=\vect{u}_h^n\) and compute
  \(\vect{w}_h^{n+1}\), \(\vect{w}_h^{n+1/2}\) from
  \eqref{eq:density-endpoint-extrapolation}
  with \(r=3\).
  \item Compute the density stages in \eqref{eq:density-ssprk3}, applying
  \(\Pi_h^{\mathrm{BP}}\) after each stage and using the limited states
  in all subsequent stage evaluations.
  \item Form \(\vect{m}_h^{n-j}=\rho_h^{n-j}\vect{u}_h^{n-j}\),
  \(j=0,1,2\), and solve
  \eqref{eq:bdf-momentum}--\eqref{eq:bdf-divergence} with \(r=3\) for
  \((\vect{u}_h^{n+1},p_h^{n+1})\), keeping \(\rho_h^{n+1}\) and
  \(\vect{w}_h^{n+1}\) fixed.
\end{enumerate}

\section{Numerical verification and benchmarks}
\label{sec:numerical-results}

The scheme is implemented in NGSolve \citep{Schoeberl2014NGSolve}.
The diagnostics include \(\|\nabla\cdot \vect{u}_h\|_{L^2(\Omega)}\),
the total mass, and the absolute mass error:
\[
  \mathcal M_h(t)=\int_\Omega\rho_h(t)\,dx,
  \qquad |\mathcal M_h(t)-\mathcal M_h(0)|.
\]
The mass denoted by \(M(t)\) on the plot axes is \(\mathcal M_h(t)\).
In the gravity-driven tests, the momentum load is
\(\vect{f}^{n+1}=\rho_h^{n+1}(0,g)^{\mathsf T}\).

\subsection{Triangular meshes: verification and benchmarks}

Unless stated otherwise, the spatial discretization uses
\(\mathrm{BDM}_2/\mathbb P_1\) for velocity and pressure and \(\mathbb P_2\)
for density (\(k_u=k_\rho=2\)).  The physical benchmarks use unstructured
triangular meshes generated by Netgen \citep{Schoeberl1997Netgen} through
NGSolve \citep{Schoeberl2014NGSolve}; the meshes used for the accuracy tests
are specified separately below.

Time integration uses SSPRK3 for density and linearly implicit BDF3 for
momentum, with stagewise BP scaling and the IPDG viscous form
\(\mathcal A_h^{\mathrm{IP}}\).  For the physical benchmarks, the initial
velocity is \(\vect{u}_h^0=\vect{0}\), and the transient Dirichlet penalty
parameter is \(\gamma_D=1\).
The nominal target time step is
\begin{equation}
  \Delta t_{\mathrm{target}}=\mathrm{CFL}_{\mathrm{nom}}\,
  \frac{h_{\mathrm{eff}}}{U_{\mathrm{ref}}},
  \label{eq:numerical-nominal-cfl}
\end{equation}
where \(\mathrm{CFL}_{\mathrm{nom}}>0\) is the nominal CFL parameter
specified for each test, \(h_{\mathrm{eff}}\) is the area-weighted
harmonic-mean mesh size, and \(U_{\mathrm{ref}}=1\) is the reference speed
unless specified otherwise.

\subsubsection{Accuracy and convergence}

To assess accuracy, we use the following variable-density exact solution
from Johnston and Liu~\citep{JohnstonLiu2002}:
\[
A(t)=0.75+0.25\sin t,\qquad
\rho=A(t)(2+\cos x\cos y),
\]
\[
\vect{u}=A(t)(-\sin^2x\sin y\cos y,\ \sin x\cos x\sin^2y),
\qquad
p=A(t)\cos x\sin y .
\]
The corresponding density source and momentum forcing are
\begin{equation}
  s_\rho
  =\partial_t\rho+\nabla\cdot(\rho\vect{u}),
  \qquad
  \vect{f}
  =\partial_t(\rho\vect{u})
   +\nabla\cdot(\rho\vect{u}\otimes\vect{u})
   -\nu\Delta\vect{u}+\nabla p .
  \label{eq:manufactured-forcing}
\end{equation}
The density residual with this source term is
\begin{equation}
  \mathcal R_h^s(q_h,\vect{w}_h,t;\varphi_h)
  :=\mathcal R_h(q_h,\vect{w}_h;\varphi_h)
    +(s_\rho(t),\varphi_h)_\Omega.
  \label{eq:manufactured-density-residual}
\end{equation}
The SSPRK3 stages evaluate \(\mathcal R_h^s\) at \(t_n\),
\(t_{n+1}\), and \(t_{n+1/2}\), respectively.  The source term is
included only in these accuracy tests; the BP estimates of
Propositions~\ref{prop:forward-euler-average-bound} and
\ref{prop:ssprk3-bp} assume \(s_\rho=0\).

The domains are \([0,2\pi]^2\) with periodic boundary conditions and
\([0,\pi]^2\) with exact Dirichlet velocity data.  In the latter case,
the exact velocity vanishes on all four walls, so
\(\Gamma_D=\partial\Omega\) and \(\vect{u}_D=\vect{0}\).  SSPRK3--BDF3 is applied up to
\(T=0.1\), with a target step size of
\(\Delta t=0.1h\).
Here \(h=2\pi/N\) in the periodic case and \(h=\pi/N\) in the
Dirichlet case.  Finite element representations of the exact
density and velocity at \(t_0,t_1,t_2\), whenever these times precede
\(T\), provide the startup history.
Tables~\ref{tab:periodic-convergence} and~\ref{tab:dirichlet-convergence}
report the errors and observed rates for the
\(\mathrm{BDM}_2/\mathbb P_1\) velocity--pressure pair with
\(\mathbb P_2\) density.  The observed rates are consistent with third-order
\(L^2\) accuracy for density and velocity and second-order accuracy for pressure.

\begin{table}[!htb]
\centering
\scriptsize
\caption{Accuracy test on \([0,2\pi]^2\) with periodic boundary conditions
and the Johnston--Liu exact solution, using SSPRK3--BDF3 with a target time
step of \(\Delta t=0.1h\) and final time \(T=0.1\).}
\resizebox{\textwidth}{!}{%
\begin{tabular}{lrrrrrrrr}
\toprule
\(\nu\) & \(N\) & \(\|\rho-\rho_h\|_{L^2}\) & order & \(\|\vect{u}-\vect{u}_h\|_{L^2}\) & order & \(\|p-p_h\|_{L^2}\) & order & \(\|\nabla\cdot \vect{u}_h\|_{L^2}\)\\
\midrule
\multirow{5}{*}{\(\pi/4000\)}
& 10  & 1.67E-02 & --   & 2.67E-02 & --   & 8.65E-02 & --   & 3.01E-15\\
& 20  & 2.77E-03 & 2.59 & 3.55E-03 & 2.91 & 1.81E-02 & 2.26 & 5.36E-15\\
& 40  & 1.74E-04 & 4.00 & 3.79E-04 & 3.23 & 3.88E-03 & 2.22 & 1.07E-14\\
& 80  & 2.01E-05 & 3.11 & 5.88E-05 & 2.69 & 9.70E-04 & 2.00 & 2.00E-14\\
& 160 & 2.43E-06 & 3.05 & 8.10E-06 & 2.86 & 2.42E-04 & 2.00 & 4.18E-14\\
\midrule
\multirow{5}{*}{1}
& 10  & 1.67E-02 & --   & 2.84E-02 & --   & 9.74E-01 & --   & 2.79E-15\\
& 20  & 2.77E-03 & 2.59 & 3.63E-03 & 2.97 & 2.88E-01 & 1.76 & 3.93E-15\\
& 40  & 1.74E-04 & 4.00 & 4.75E-04 & 2.94 & 7.52E-02 & 1.94 & 9.28E-15\\
& 80  & 2.01E-05 & 3.11 & 5.95E-05 & 3.00 & 1.93E-02 & 1.97 & 1.88E-14\\
& 160 & 2.43E-06 & 3.05 & 7.43E-06 & 3.00 & 4.84E-03 & 1.99 & 6.41E-13\\
\bottomrule
\end{tabular}}
\label{tab:periodic-convergence}
\end{table}

\begin{table}[!htb]
\centering
\scriptsize
\caption{Accuracy test on \([0,\pi]^2\) with Dirichlet boundary conditions
and the Johnston--Liu exact solution, using SSPRK3--BDF3 with a target time
step of \(\Delta t=0.1h\) and final time \(T=0.1\).}
\resizebox{\textwidth}{!}{%
\begin{tabular}{lrrrrrrrr}
\toprule
\(\nu\) & \(N\) & \(\|\rho-\rho_h\|_{L^2}\) & order & \(\|\vect{u}-\vect{u}_h\|_{L^2}\) & order & \(\|p-p_h\|_{L^2}\) & order & \(\|\nabla\cdot \vect{u}_h\|_{L^2}\)\\
\midrule
\multirow{5}{*}{\(\pi/4000\)}
& 10  & 1.36E-03 & --   & 1.78E-03 & --   & 1.05E-02 & --   & 3.33E-15\\
& 20  & 8.65E-05 & 3.97 & 2.06E-04 & 3.11 & 1.94E-03 & 2.44 & 6.17E-15\\
& 40  & 1.01E-05 & 3.10 & 3.05E-05 & 2.76 & 4.85E-04 & 2.00 & 1.23E-14\\
& 80  & 1.25E-06 & 3.02 & 4.12E-06 & 2.89 & 1.21E-04 & 2.00 & 2.61E-14\\
& 160 & 1.43E-07 & 3.13 & 4.95E-07 & 3.06 & 3.03E-05 & 2.00 & 5.33E-14\\
\midrule
\multirow{5}{*}{1}
& 10  & 1.36E-03 & --   & 1.86E-03 & --   & 1.50E-01 & --   & 1.64E-15\\
& 20  & 8.63E-05 & 3.98 & 2.44E-04 & 2.93 & 3.81E-02 & 1.98 & 3.80E-15\\
& 40  & 1.01E-05 & 3.09 & 3.00E-05 & 3.02 & 9.68E-03 & 1.98 & 9.07E-15\\
& 80  & 1.25E-06 & 3.02 & 3.72E-06 & 3.01 & 2.43E-03 & 2.00 & 1.64E-14\\
& 160 & 1.43E-07 & 3.12 & 4.64E-07 & 3.00 & 6.07E-04 & 2.00 & 3.51E-14\\
\bottomrule
\end{tabular}}
\label{tab:dirichlet-convergence}
\end{table}

\subsubsection{Low-density bound-preservation test}
\label{subsubsec:lowdensity-bdm}

To assess the limiter in the low-density regime, we use the exact solution of
Zhang, Xia, and Shu~\citep[Example~6.6]{ZhangXiaShu2012JSC} on the doubly periodic square
\(\Omega=[0,1]^2\):
\[
  \rho=1+0.99\sin\bigl(2\pi(x+y-2t)\bigr),\qquad
  \vect{u}=(1,1)^{\mathsf T},\qquad p=0 ,
\]
so that the density oscillates between \(0.01\) and \(1.99\).  The exact density
is smooth and the data satisfy the homogeneous density and momentum equations
with \(s_\rho=0\) and \(\vect{f}=\vect{0}\), so the test isolates the effect of
the limiter.  The viscosity is \(\nu=\pi/4000\).

The density is limited at every SSPRK3 stage with the triangular control set of
Zhang, Xia, and Shu~\citep[Section~3]{ZhangXiaShu2012JSC}, the scaling \eqref{eq:bp-scaling}, and the
bounds \(m=\delta\), \(M=2-\delta\) with \(\delta=0.01\).  SSPRK3--BDF3 is
applied up to \(T=0.1\) with a target step size of \(\Delta t=0.01/N\) on
unstructured triangular meshes generated by Netgen \citep{Schoeberl1997Netgen}
through NGSolve \citep{Schoeberl2014NGSolve}.  Table~\ref{tab:lowdensity-bdm} reports the errors and observed rates.

\begin{table}[!htb]
\centering
\scriptsize
\caption{Low-density test on the unstructured triangular meshes}
\resizebox{\textwidth}{!}{%
\begin{tabular}{llrrrrrrr}
\toprule
startup & \(N\) & \(\|\rho-\rho_h\|_{L^1}\) & order & \(\|\rho-\rho_h\|_{L^2}\) & order & \(\|\rho-\rho_h\|_{L^\infty}\) & order & \(\|\nabla\cdot\vect{u}_h\|_{L^2}\)\\
\midrule
\multirow{4}{*}{exact history}
 & 10  & 1.70E-3 & --   & 3.22E-3 & --   & 6.41E-2 & --   & 1.91E-13\\
 & 20  & 1.01E-4 & 4.05 & 1.47E-4 & 4.43 & 2.66E-3 & 4.56 & 3.81E-13\\
 & 40  & 1.28E-5 & 3.01 & 1.97E-5 & 2.93 & 4.70E-4 & 2.53 & 4.86E-11\\
 & 60  & 3.48E-6 & 3.17 & 5.34E-6 & 3.18 & 1.37E-4 & 3.00 & 4.39E-11\\
\midrule
\multirow{4}{*}{SSPRK3 startup}
 & 10  & 1.69E-3 & --   & 3.21E-3 & --   & 6.41E-2 & --   & 1.92E-13\\
 & 20  & 1.01E-4 & 4.04 & 1.47E-4 & 4.42 & 2.66E-3 & 4.56 & 3.82E-13\\
 & 40  & 1.28E-5 & 3.01 & 1.97E-5 & 2.93 & 4.70E-4 & 2.53 & 7.47E-13\\
 & 60  & 3.49E-6 & 3.17 & 5.34E-6 & 3.18 & 1.37E-4 & 3.00 & 2.68E-11\\
\bottomrule
\end{tabular}}
\label{tab:lowdensity-bdm}
\end{table}

\FloatBarrier
\subsubsection{Rayleigh--Taylor}

The Rayleigh--Taylor configuration follows Tryggvason~\citep{Tryggvason1988RT}
and Li et al.~\citep{Li2021BP}
on
\(\Omega=[-0.5,0.5]\times[-2,1]\), with \(\vect{u}_0=\vect{0}\).
All four walls are no-slip, with \(\Gamma_D=\partial\Omega\) and
\(\vect{u}_D=\vect{0}\).  The smoothed initial density is
\begin{equation}
  \eta(x)=-a\cos\left(\frac{2\pi x}{d}\right),
  \qquad
  \rho_0(x,y)=
  \frac{\rho_M+\rho_m}{2}
  +
  \frac{\rho_M-\rho_m}{2}
  \tanh\left(\frac{y-\eta(x)}{\varepsilon}\right),
  \label{eq:rt-initial-density}
\end{equation}
with \(a=0.1\), \(d=1\), \(\rho_m=1\), and
\(\varepsilon=0.05\).  For the Atwood number
\[
  A_t=\frac{\rho_M-\rho_m}{\rho_M+\rho_m},
  \qquad
  \rho_M=\rho_m\,\frac{1+A_t}{1-A_t},
\]
the values \(A_t=0.5,0.75,0.9\) give the admissible density
intervals \([1,3]\), \([1,7]\), and \([1,19]\).

\begin{figure}[!htb]
\centering
\begin{minipage}{0.32\textwidth}
\centering
\includegraphics[width=\linewidth]{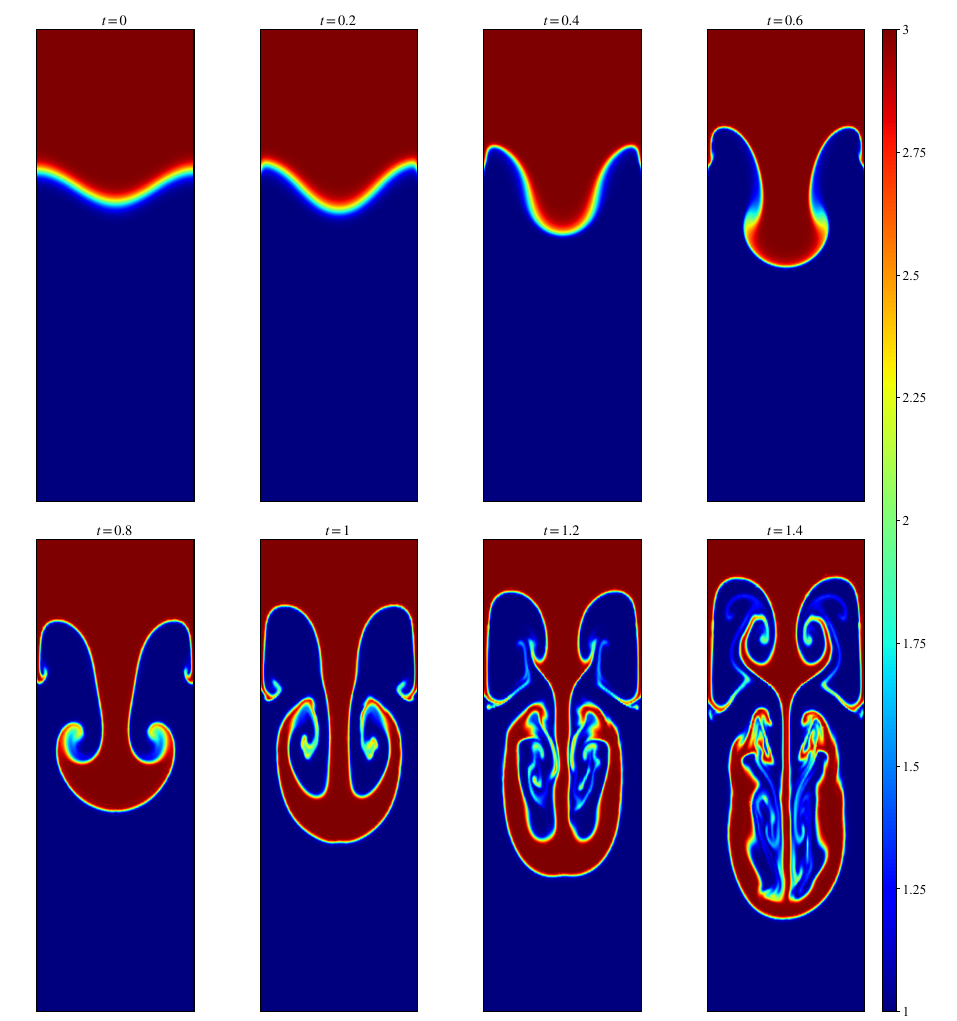}\par
\small (a) \(A_t=0.5\)
\end{minipage}\hfill
\begin{minipage}{0.32\textwidth}
\centering
\includegraphics[width=\linewidth]{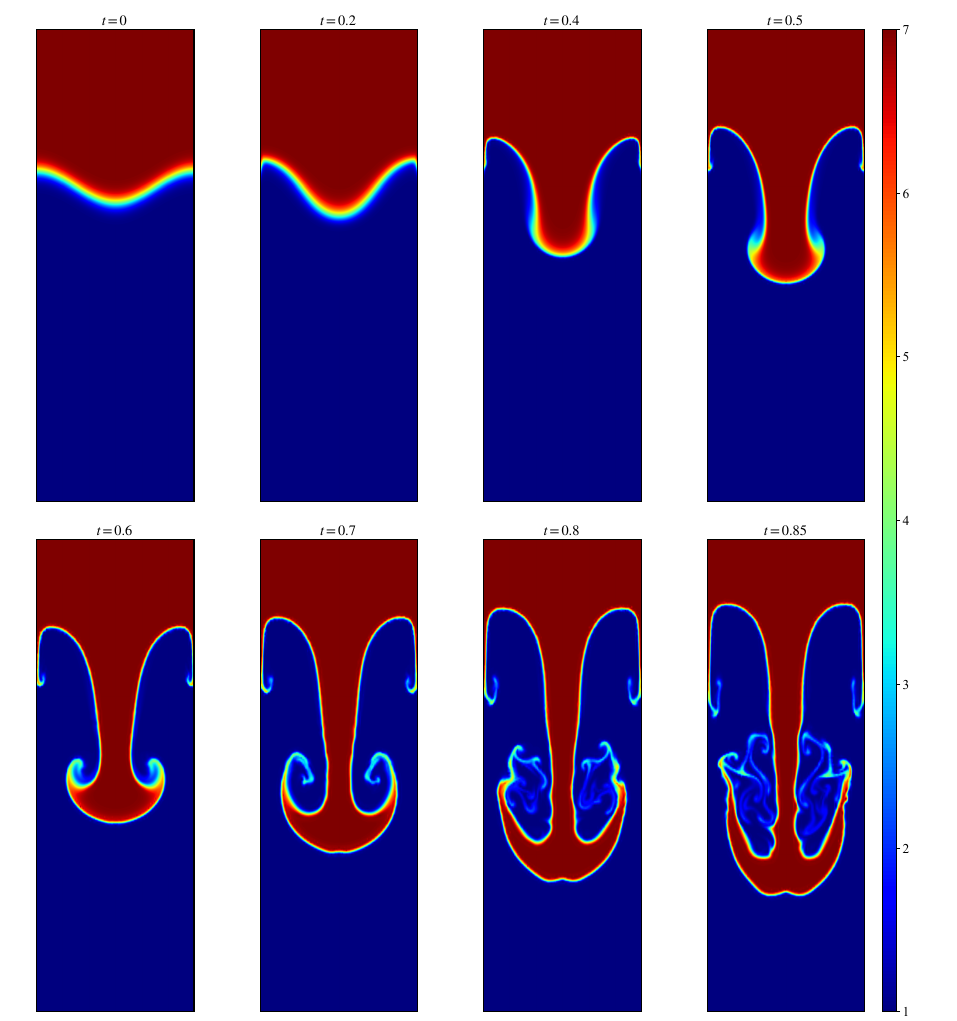}\par
\small (b) \(A_t=0.75\)
\end{minipage}\hfill
\begin{minipage}{0.32\textwidth}
\centering
\includegraphics[width=\linewidth]{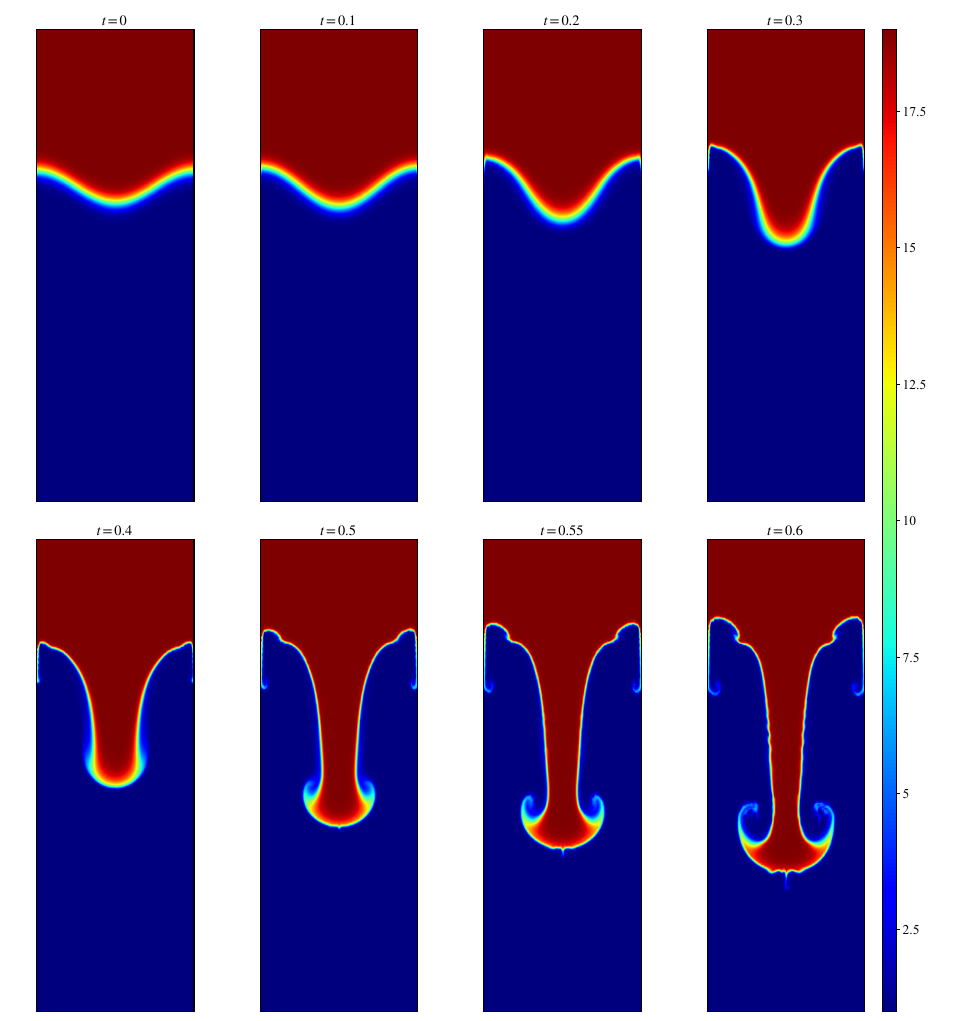}\par
\small (c) \(A_t=0.9\)
\end{minipage}
\caption{Rayleigh--Taylor density evolution for the three Atwood numbers
\(A_t=0.5,0.75,0.9\), shown from left to right at the times specified in
the corresponding simulations.}
\label{fig:rt-density-combined}
\end{figure}

\begin{figure}[!htb]
\centering
\begin{minipage}{0.32\textwidth}
\centering
\includegraphics[width=\textwidth]{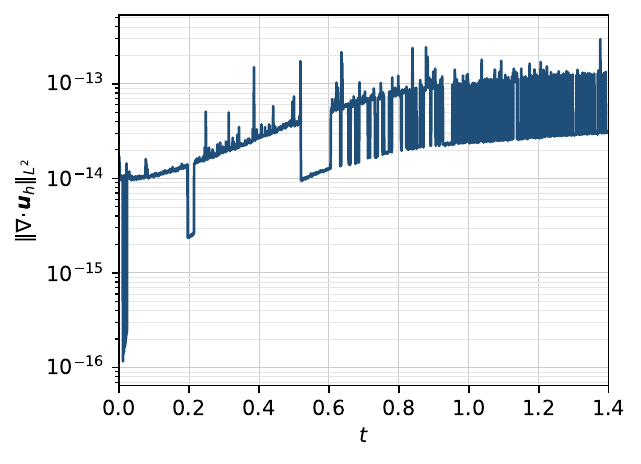}\\
\small \(A_t=0.5\)
\end{minipage}\hfill
\begin{minipage}{0.32\textwidth}
\centering
\includegraphics[width=\textwidth]{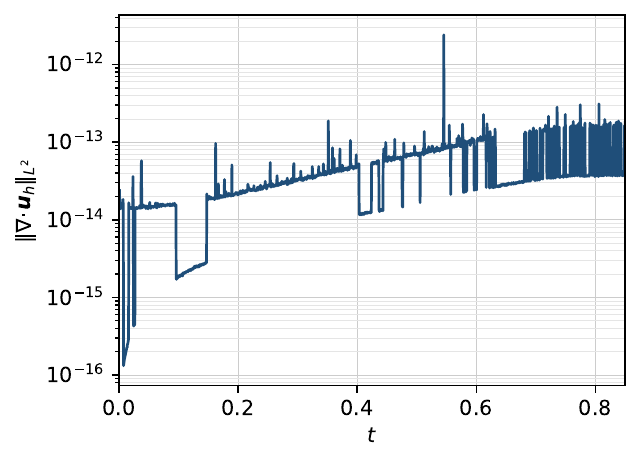}\\
\small \(A_t=0.75\)
\end{minipage}\hfill
\begin{minipage}{0.32\textwidth}
\centering
\includegraphics[width=\textwidth]{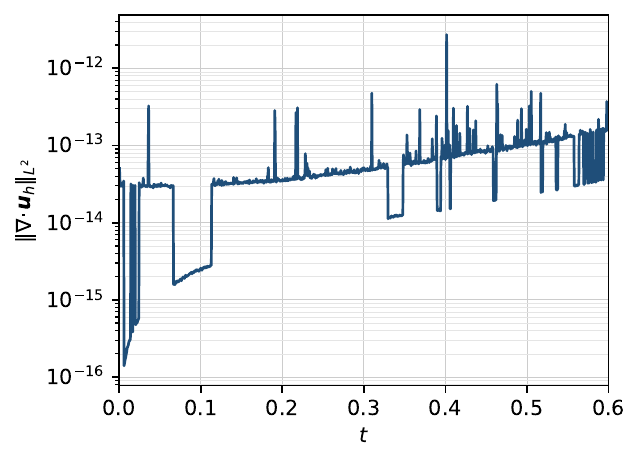}\\
\small \(A_t=0.9\)
\end{minipage}
\caption{Divergence diagnostics for the displayed Rayleigh--Taylor runs on a
base-10 logarithmic scale, ordered as \(A_t=0.5,0.75,0.9\).}
\label{fig:rt-div-monitor-combined}
\end{figure}

\begin{figure}[!htb]
\centering
\begin{minipage}{0.32\textwidth}
\centering
\includegraphics[width=\textwidth]{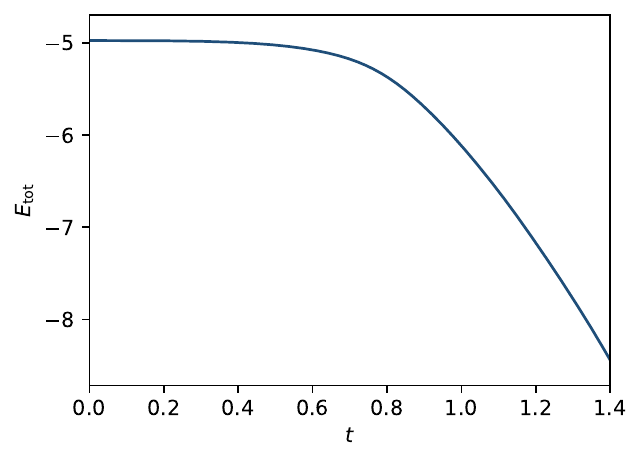}\\
\small \(A_t=0.5\)
\end{minipage}\hfill
\begin{minipage}{0.32\textwidth}
\centering
\includegraphics[width=\textwidth]{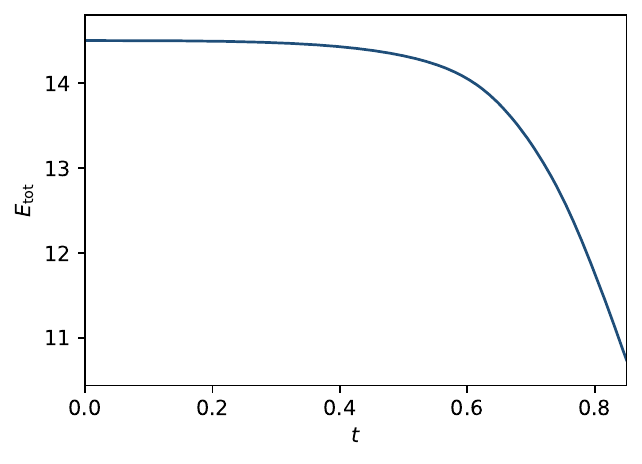}\\
\small \(A_t=0.75\)
\end{minipage}\hfill
\begin{minipage}{0.32\textwidth}
\centering
\includegraphics[width=\textwidth]{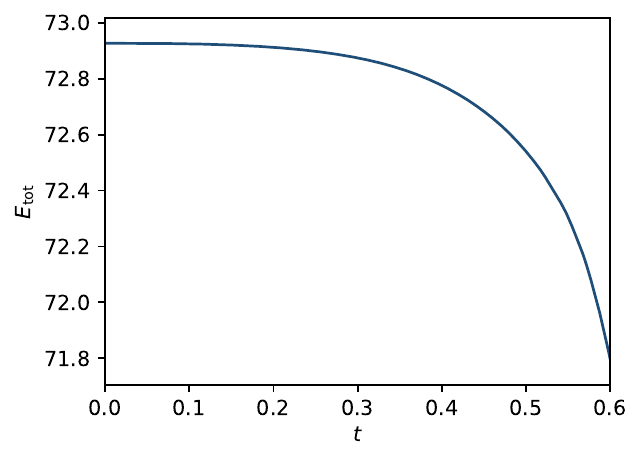}\\
\small \(A_t=0.9\)
\end{minipage}
\caption{Total-energy histories \(E_{\mathrm{tot}}\) for the displayed
Rayleigh--Taylor runs, ordered as \(A_t=0.5,0.75,0.9\); see
\eqref{eq:rt-energy-definition}.}
\label{fig:rt-energy-total-combined}
\end{figure}

\begin{figure}[!htb]
\centering
\begin{minipage}{0.32\textwidth}
\centering
\includegraphics[width=\textwidth]{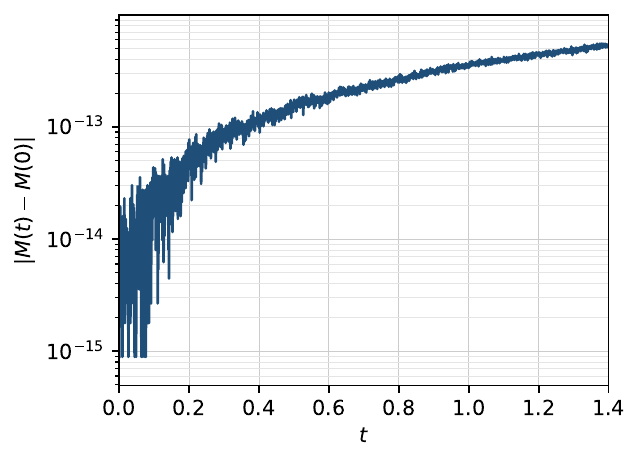}\\
\small \(A_t=0.5\)
\end{minipage}\hfill
\begin{minipage}{0.32\textwidth}
\centering
\includegraphics[width=\textwidth]{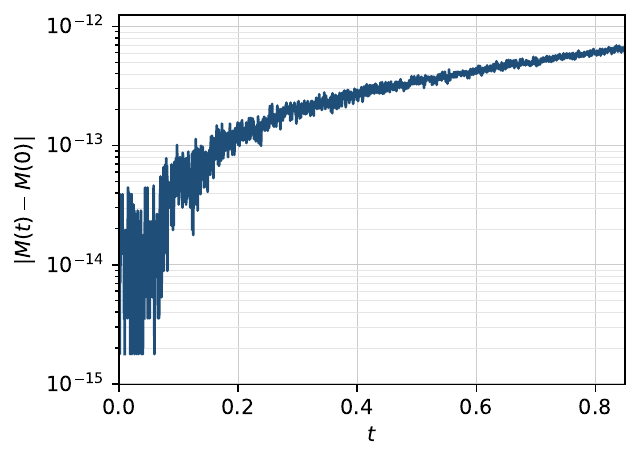}\\
\small \(A_t=0.75\)
\end{minipage}\hfill
\begin{minipage}{0.32\textwidth}
\centering
\includegraphics[width=\textwidth]{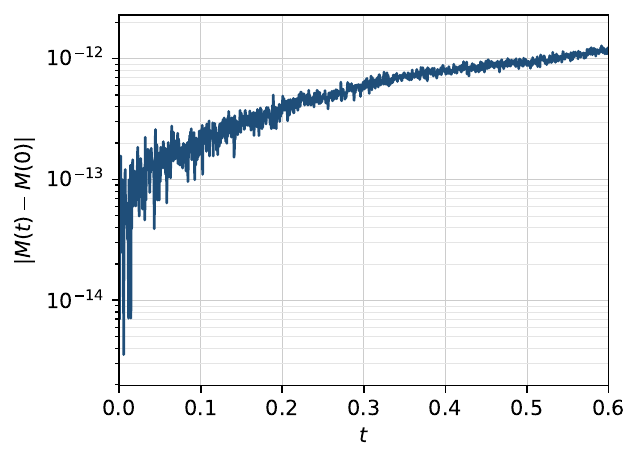}\\
\small \(A_t=0.9\)
\end{minipage}
\caption{Absolute mass-conservation diagnostics for the displayed
Rayleigh--Taylor runs on a base-10 logarithmic scale, ordered as
\(A_t=0.5,0.75,0.9\).}
\label{fig:rt-mass-error-combined}
\end{figure}

\begin{figure}[!htb]
\centering
\begin{minipage}{0.32\textwidth}
\centering
\includegraphics[width=\linewidth]{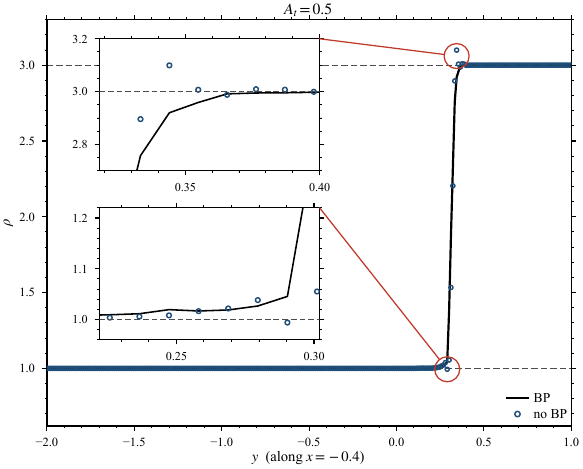}\par
\small (a) \(A_t=0.5\)
\end{minipage}\hfill
\begin{minipage}{0.32\textwidth}
\centering
\includegraphics[width=\linewidth]{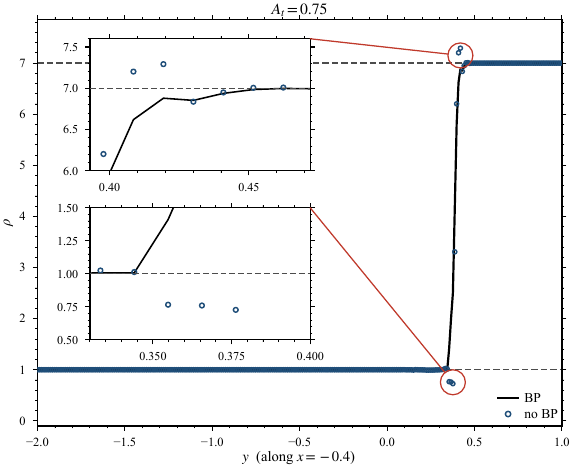}\par
\small (b) \(A_t=0.75\)
\end{minipage}\hfill
\begin{minipage}{0.32\textwidth}
\centering
\includegraphics[width=\linewidth]{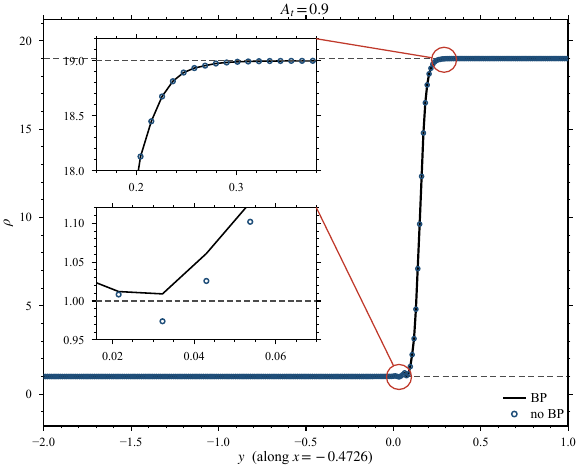}\par
\small (c) \(A_t=0.9\)
\end{minipage}
\caption{Density profiles for the Rayleigh--Taylor runs with
\(A_t=0.5,0.75,0.9\), taken at \(t=0.5,0.5,0.15\), respectively.
The solid lines denote the solutions with BP limiting and the open circles
the solutions without BP limiting; the dashed lines mark the prescribed
density bounds.  The first two profiles are taken along \(x=-0.4\), while
the third is taken along \(x=-0.4726\).}
\label{fig:rt-bp-profiles-combined}
\end{figure}

All runs use \(g=-9.80665\), \(Re=1000\), and
\[
  \nu=\frac{\rho_m d^{3/2}\sqrt{|g|}}{Re}.
\]
The nominal mesh spacing is \(h=1/40\), corresponding to
\(N_x=40\), \(N_y=120\) across the domain; the interior triangulation
is unstructured.  We set \(\mathrm{CFL}_{\mathrm{nom}}=0.02\).
The energy diagnostics are
\begin{equation}
  E_{\mathrm{kin}}(t)
  =
  \frac12\int_\Omega \rho_h |\vect{u}_h|^2\,dx,
  \qquad
  E_{\mathrm{pot}}(t)
  =
  \int_\Omega \rho_h(-g)y\,dx,
  \qquad
  E_{\mathrm{tot}}=E_{\mathrm{kin}}+E_{\mathrm{pot}}.
  \label{eq:rt-energy-definition}
\end{equation}
The potential energy is measured relative to \(y=0\).
Figure~\ref{fig:rt-density-combined} shows the density evolution, while
Figures~\ref{fig:rt-div-monitor-combined}--\ref{fig:rt-mass-error-combined}
report the divergence, energy, and mass diagnostics.
Figure~\ref{fig:rt-bp-profiles-combined} compares the limited and unlimited
density profiles and shows the prescribed density bounds.

\subsubsection{Rising bubble}

The rising-bubble configuration is taken from Li et al.~\citep{Li2021BP}
and Li et al.~\citep{Li2024Energy}.
For this test, the constant-viscosity term \(-\nu\Delta\vect{u}\) is
replaced by \(-\nabla\cdot(\mu(\rho)\nabla\vect{u})\), as in the
variable-viscosity extension in Section~4.2 of Li et al.~\citep{Li2024Energy}.

The domain is \(\Omega=[-0.01,0.01]\times[0,0.03]\), with
\(\vect{u}_0=\vect{0}\) and no-slip conditions on all four walls.
Thus \(\Gamma_D=\partial\Omega\) and \(\vect{u}_D=\vect{0}\).  Define
\[
  r(x,y)=\sqrt{(x-c_x)^2+(y-c_y)^2}.
\]
The smoothed initial density is
\begin{equation}
  \rho_0(x,y)
  =
  \frac{\rho_w+\rho_a}{2}
  +
  \frac{\rho_w-\rho_a}{2}
  \tanh\left(\frac{r(x,y)-R_0}{\varepsilon}\right),
  \label{eq:rising-bubble-density}
\end{equation}
with \(\rho_a=1.161\), \(\rho_w=995.65\),
\((c_x,c_y)=(0,0.0075)\), \(R_0=0.0025\), and
\(\varepsilon=h=2.50\times10^{-4}\).  The viscosity depends affinely on
the density:
\begin{equation}
  \mu(\rho_h)
  =
  \mu_a
  +
  \frac{\mu_w-\mu_a}{\rho_w-\rho_a}(\rho_h-\rho_a),
  \qquad
  \mu_a=1.86\times10^{-5},\quad
  \mu_w=7.977\times10^{-4}.
  \label{eq:rising-bubble-viscosity}
\end{equation}

\begin{figure}[!htb]
\centering
\begin{minipage}{0.33\textwidth}
\centering
\includegraphics[width=\textwidth]{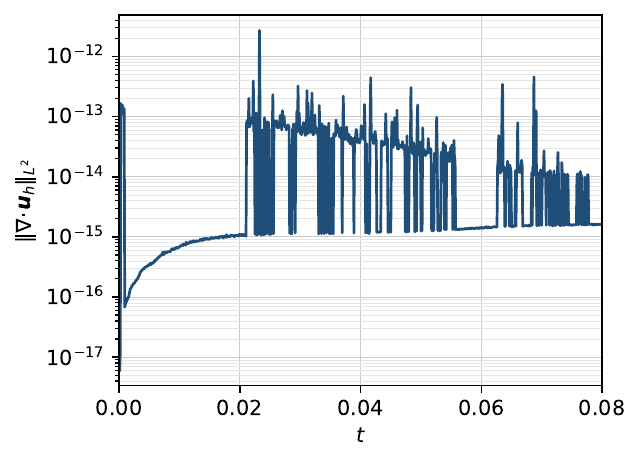}\\
\small (a) divergence
\end{minipage}\hfill
\begin{minipage}{0.33\textwidth}
\centering
\includegraphics[width=\textwidth]{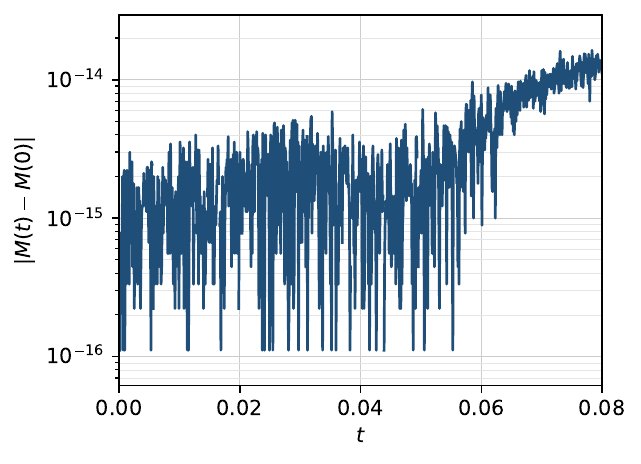}\\
\small (b) absolute mass error
\end{minipage}\hfill
\begin{minipage}{0.33\textwidth}
\centering
\includegraphics[width=\textwidth]{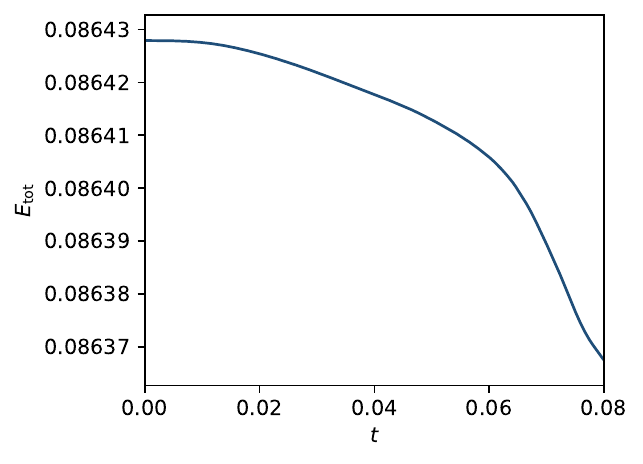}\\
\small (c) total energy
\end{minipage}
\caption{Rising-bubble diagnostics:
(a) \(\|\nabla\cdot \vect{u}_h\|_{L^2}\), (b) the absolute mass error
\(|\mathcal M_h(t)-\mathcal M_h(0)|\), and (c) total energy \(E_{\mathrm{tot}}\).  The first two
panels use base-10 logarithmic vertical axes.}
\label{fig:rising-diagnostics}
\end{figure}

\begin{figure}[!htb]
\centering
\includegraphics[width=0.5\textwidth]{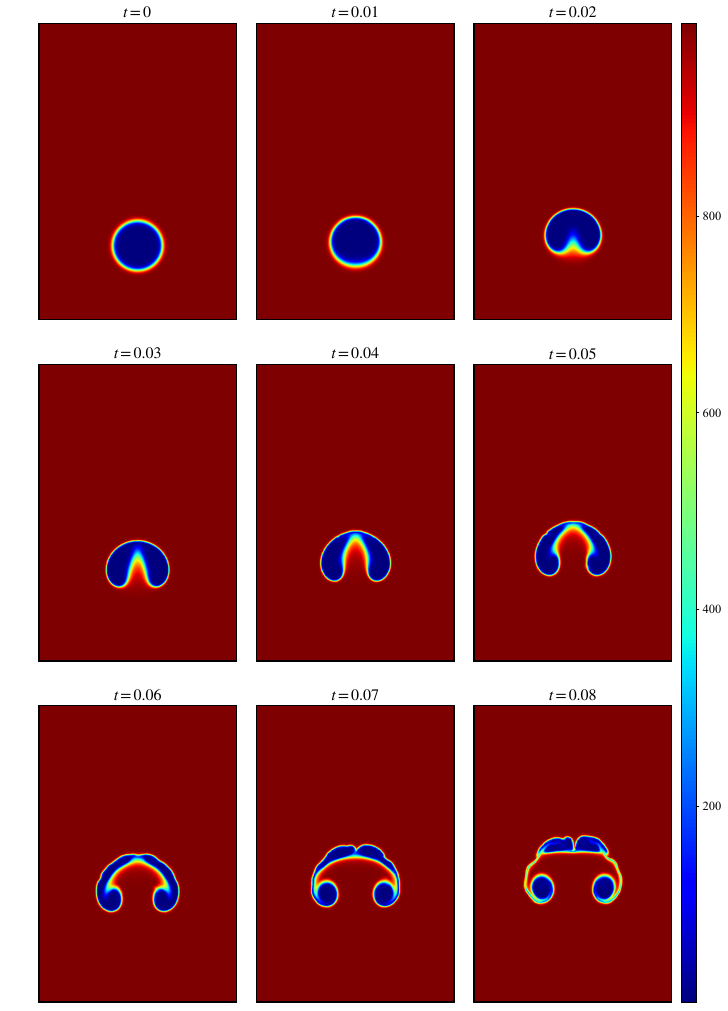}
\caption{Rising-bubble density evolution at different times.}
\label{fig:rising-density-panel}
\end{figure}

The calculation uses \(g=-9.80665\), \(h_{\max}=2.50\times10^{-4}\), and \(\mathrm{CFL}_{\mathrm{nom}}=0.4\).
Figure~\ref{fig:rising-density-panel} shows the ascent and deformation of
the light bubble.  The divergence, mass, and energy histories appear in
Figure~\ref{fig:rising-diagnostics}.  The recorded density averages
remain in \([1.161,995.65]\) up to roundoff.

\subsection{Rectangular meshes: verification and benchmarks}
\label{subsec:rt-meshes}
\label{subsec:rt-slip-lock-exchange}

The computations in this subsection use the \(\mathrm{RT}_2/\mathbb Q_2\)
velocity--pressure pair with discontinuous \(\mathbb Q_2\) density on
structured rectangular meshes, together with the SSPRK3--BDF3 time
integration and the BP limiter described in
Sections~\ref{subsec:density-equation}--\ref{subsec:complete-algorithm}.

\subsubsection{Accuracy and convergence}
\label{subsubsec:accuracy-rt}

The accuracy test of the preceding subsection is repeated with the
\(\mathrm{RT}_2/\mathbb Q_2\) velocity--pressure pair and a discontinuous
\(\mathbb Q_2\) density on structured \(N\times N\) quadrilateral meshes, so
that the mesh consists of \(N^2\) rectangles of side \(h=2\pi/N\) in the
periodic case and \(h=\pi/N\) in the Dirichlet case.  The exact solution,
viscosities, source terms, final time, target step size and startup history
are those of the triangular accuracy test.  On rectangles the limiter uses
the tensor-product control set, positive quadrature, and corresponding face
CFL condition described by Zhang and Shu~\citep[Section~3.2]{ZhangShu2010}.
Tables~\ref{tab:periodic-convergence-rt} and~\ref{tab:dirichlet-convergence-rt}
report the errors and rates.  
\begin{table}[!htb]
\centering
\scriptsize
\caption{Accuracy test on \([0,2\pi]^2\) with periodic boundary conditions
and the Johnston--Liu exact solution for the \(\mathrm{RT}_2/\mathbb Q_2\)
discretisation on structured \(N\times N\) quadrilateral meshes, using
SSPRK3--BDF3 with a target time step of \(\Delta t=0.1h\) and final time \(T=0.1\).}
\resizebox{\textwidth}{!}{%
\begin{tabular}{lrrrrrrrr}
\toprule
\(\nu\) & \(N\) & \(\|\rho-\rho_h\|_{L^2}\) & order & \(\|\vect{u}-\vect{u}_h\|_{L^2}\) & order & \(\|p-p_h\|_{L^2}\) & order & \(\|\nabla\cdot \vect{u}_h\|_{L^2}\)\\
\midrule
\multirow{5}{*}{\(\pi/4000\)}
& 10  & 2.71E-03 & --   & 9.21E-03 & --   & 2.75E-03 & --   & 1.95E-16\\
& 20  & 3.48E-04 & 2.96 & 1.18E-03 & 2.96 & 3.39E-04 & 3.02 & 5.51E-16\\
& 40  & 4.94E-05 & 2.82 & 1.72E-04 & 2.78 & 4.23E-05 & 3.00 & 9.59E-16\\
& 80  & 7.59E-06 & 2.70 & 2.85E-05 & 2.59 & 5.28E-06 & 3.00 & 2.08E-15\\
& 160 & 1.02E-06 & 2.90 & 3.94E-06 & 2.86 & 6.59E-07 & 3.00 & 6.20E-15\\
\midrule
\multirow{5}{*}{1}
& 10  & 2.71E-03 & --   & 1.29E-02 & --   & 1.53E-02 & --   & 1.88E-16\\
& 20  & 3.50E-04 & 2.95 & 1.84E-03 & 2.81 & 1.55E-03 & 3.30 & 5.31E-16\\
& 40  & 4.97E-05 & 2.82 & 2.27E-04 & 3.02 & 2.06E-04 & 2.91 & 9.61E-16\\
& 80  & 7.62E-06 & 2.70 & 2.84E-05 & 3.00 & 2.57E-05 & 3.00 & 2.13E-15\\
& 160 & 1.02E-06 & 2.90 & 3.55E-06 & 3.00 & 3.21E-06 & 3.00 & 6.04E-15\\
\bottomrule
\end{tabular}}
\label{tab:periodic-convergence-rt}
\end{table}

\begin{table}[!htb]
\centering
\scriptsize
\caption{Accuracy test on \([0,\pi]^2\) with Dirichlet boundary conditions
and the Johnston--Liu exact solution for the \(\mathrm{RT}_2/\mathbb Q_2\)
discretisation on structured \(N\times N\) quadrilateral meshes, using
SSPRK3--BDF3 with a target time step of \(\Delta t=0.1h\) and final time \(T=0.1\).}
\resizebox{\textwidth}{!}{%
\begin{tabular}{lrrrrrrrr}
\toprule
\(\nu\) & \(N\) & \(\|\rho-\rho_h\|_{L^2}\) & order & \(\|\vect{u}-\vect{u}_h\|_{L^2}\) & order & \(\|p-p_h\|_{L^2}\) & order & \(\|\nabla\cdot \vect{u}_h\|_{L^2}\)\\
\midrule
\multirow{5}{*}{\(\pi/4000\)}
& 10  & 1.74E-04 & --   & 6.77E-04 & --   & 3.97E-04 & --   & 3.57E-16\\
& 20  & 2.47E-05 & 2.82 & 9.05E-05 & 2.90 & 2.17E-05 & 4.19 & 5.27E-16\\
& 40  & 3.80E-06 & 2.70 & 1.45E-05 & 2.65 & 2.64E-06 & 3.04 & 1.49E-15\\
& 80  & 5.09E-07 & 2.90 & 1.98E-06 & 2.87 & 3.30E-07 & 3.00 & 2.21E-15\\
& 160 & 6.37E-08 & 3.00 & 2.39E-07 & 3.05 & 4.12E-08 & 3.00 & 6.42E-15\\
\midrule
\multirow{5}{*}{1}
& 10  & 1.75E-04 & --   & 9.05E-04 & --   & 8.25E-04 & --   & 3.33E-16\\
& 20  & 2.48E-05 & 2.82 & 1.12E-04 & 3.01 & 1.02E-04 & 3.01 & 5.74E-16\\
& 40  & 3.81E-06 & 2.70 & 1.41E-05 & 2.99 & 1.27E-05 & 3.00 & 1.68E-15\\
& 80  & 5.10E-07 & 2.90 & 1.77E-06 & 2.99 & 1.60E-06 & 3.00 & 2.07E-15\\
& 160 & 6.37E-08 & 3.00 & 2.21E-07 & 3.00 & 2.00E-07 & 3.00 & 8.95E-15\\
\bottomrule
\end{tabular}}
\label{tab:dirichlet-convergence-rt}
\end{table}

\FloatBarrier
\subsubsection{Low-density bound-preservation test}
\label{subsubsec:lowdensity-rt}

The low-density test of Section~\ref{subsubsec:lowdensity-bdm} is repeated with
the \(\mathrm{RT}_2/\mathbb Q_2\) velocity--pressure pair and a discontinuous
\(\mathbb Q_2\) density on structured \(N\times N\) quadrilateral meshes, with
the same exact solution, viscosity, final time, target step size and startup
history.  On rectangles the limiter uses the tensor-product control set,
positive quadrature, and corresponding face CFL condition described in
Zhang and Shu~\citep[Section~3.2]{ZhangShu2010}.
Table~\ref{tab:lowdensity-rt} reports the errors and observed rates.

\begin{table}[!htb]
\centering
\scriptsize
\caption{Low-density  test of the \(\mathrm{RT}_2/\mathbb Q_2\)
discretisation on structured \(N\times N\) quadrilateral meshes}
\resizebox{\textwidth}{!}{%
\begin{tabular}{llrrrrrrr}
\toprule
startup & \(N\) & \(\|\rho-\rho_h\|_{L^1}\) & order & \(\|\rho-\rho_h\|_{L^2}\) & order & \(\|\rho-\rho_h\|_{L^\infty}\) & order & \(\|\nabla\cdot\vect{u}_h\|_{L^2}\)\\
\midrule
\multirow{4}{*}{exact history}
 & 10  & 1.21E-3 & --   & 1.58E-3 & --   & 9.06E-3 & --   & 2.54E-15\\
 & 20  & 1.14E-4 & 3.41 & 1.53E-4 & 3.36 & 1.01E-3 & 3.16 & 5.65E-15\\
 & 40  & 1.37E-5 & 3.06 & 1.88E-5 & 3.02 & 1.28E-4 & 2.99 & 1.15E-14\\
 & 60  & 4.03E-6 & 3.01 & 5.57E-6 & 3.00 & 3.79E-5 & 3.00 & 2.80E-14\\
\midrule
\multirow{4}{*}{SSPRK3 startup}
 & 10  & 1.17E-3 & --   & 1.52E-3 & --   & 9.35E-3 & --   & 2.93E-15\\
 & 20  & 1.15E-4 & 3.36 & 1.54E-4 & 3.30 & 1.01E-3 & 3.21 & 5.62E-15\\
 & 40  & 1.37E-5 & 3.06 & 1.88E-5 & 3.03 & 1.28E-4 & 2.99 & 1.27E-14\\
 & 60  & 4.03E-6 & 3.01 & 5.57E-6 & 3.01 & 3.79E-5 & 3.00 & 2.48E-14\\
\bottomrule
\end{tabular}}
\label{tab:lowdensity-rt}
\end{table}

\FloatBarrier
\subsubsection{Falling bubble}
\label{subsubsec:falling-bubble-rt}

The falling-bubble test follows Li et al.~\citep{Li2021BP}
and Li et al.~\citep{Li2024Energy} on
\(\Omega=[0,1]\times[0,2]\), with
\(\vect{u}_0=\vect{0}\) and no-slip conditions on all four walls.
Thus \(\Gamma_D=\partial\Omega\) and \(\vect{u}_D=\vect{0}\).  Define
\[
  H_\varepsilon(s)=\frac12\left(1+\tanh\frac{s}{\varepsilon}\right),
  \qquad
  r(x,y)=\sqrt{(x-c_x)^2+(y-c_y)^2}.
\]
The initial density is
\begin{equation}
  \rho_{\mathrm{layer}}
  =
  \rho_L+(\rho_H-\rho_L)H_\varepsilon(y_0-y),
  \qquad
  \rho_{\mathrm{drop}}
  =
  \rho_L+(\rho_H-\rho_L)H_\varepsilon(r_0-r),
  \qquad
  \rho_0=\max\{\rho_{\mathrm{layer}},\rho_{\mathrm{drop}}\}.
  \label{eq:falling-bubble-density}
\end{equation}
Here \(\rho_H=100\), \(\rho_L=1\), \(y_0=1\),
\((c_x,c_y)=(0.5,1.75)\), \(r_0=0.2\), and \(\varepsilon=h\).

The parameters are \(g=-1\), \(Re=1000\), and
\[
  \nu=\frac{\rho_L d^{3/2}\sqrt{|g|}}{Re},
  \qquad d=2.
\]
The discretization is the \(\mathrm{RT}_2/\mathbb Q_2\)
velocity--pressure pair with discontinuous \(\mathbb Q_2\) density on a
structured \(50\times100\) rectangular mesh, so that \(h=1/50\), and we
set \(\mathrm{CFL}_{\mathrm{nom}}=0.05\), giving \(\Delta t=10^{-3}\) up to
\(T=1.3\).  The density bounds are enforced by the BP limiter of
Section~\ref{subsubsec:density-bp}.

\begin{figure}[!htb]
\centering
\begin{minipage}{0.33\textwidth}
\centering
\includegraphics[width=\textwidth]{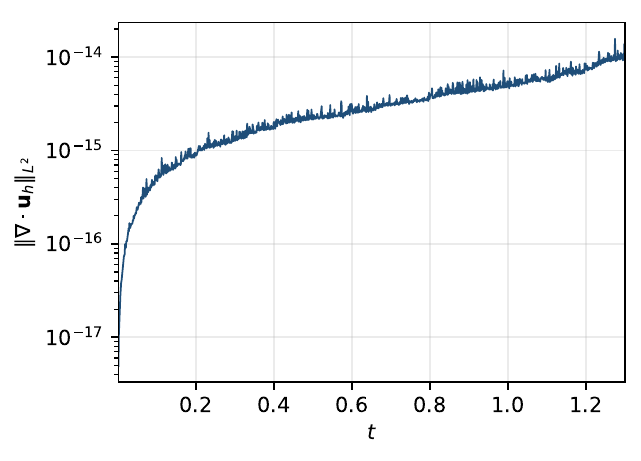}\\
\small (a) divergence
\end{minipage}\hfill
\begin{minipage}{0.33\textwidth}
\centering
\includegraphics[width=\textwidth]{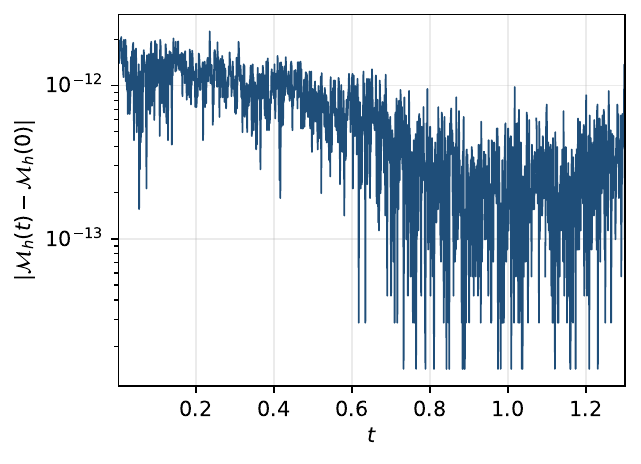}\\
\small (b) absolute mass error
\end{minipage}\hfill
\begin{minipage}{0.33\textwidth}
\centering
\includegraphics[width=\textwidth]{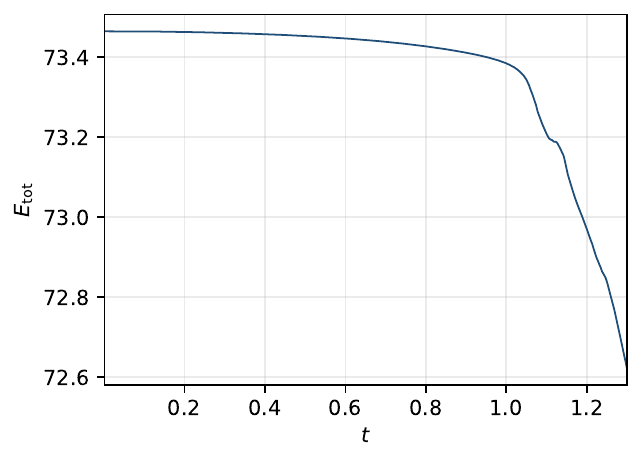}\\
\small (c) total energy
\end{minipage}
\caption{Falling-bubble diagnostics on the structured \(50\times100\)
RT mesh:
(a) \(\|\nabla\cdot \vect{u}_h\|_{L^2}\), (b) the absolute mass error
\(|\mathcal M_h(t)-\mathcal M_h(0)|\), and (c) total energy \(E_{\mathrm{tot}}\).  The first two
panels use base-10 logarithmic vertical axes.}
\label{fig:falling-diagnostics}
\end{figure}

\begin{figure}[!htb]
\centering
\includegraphics[width=0.5\textwidth]{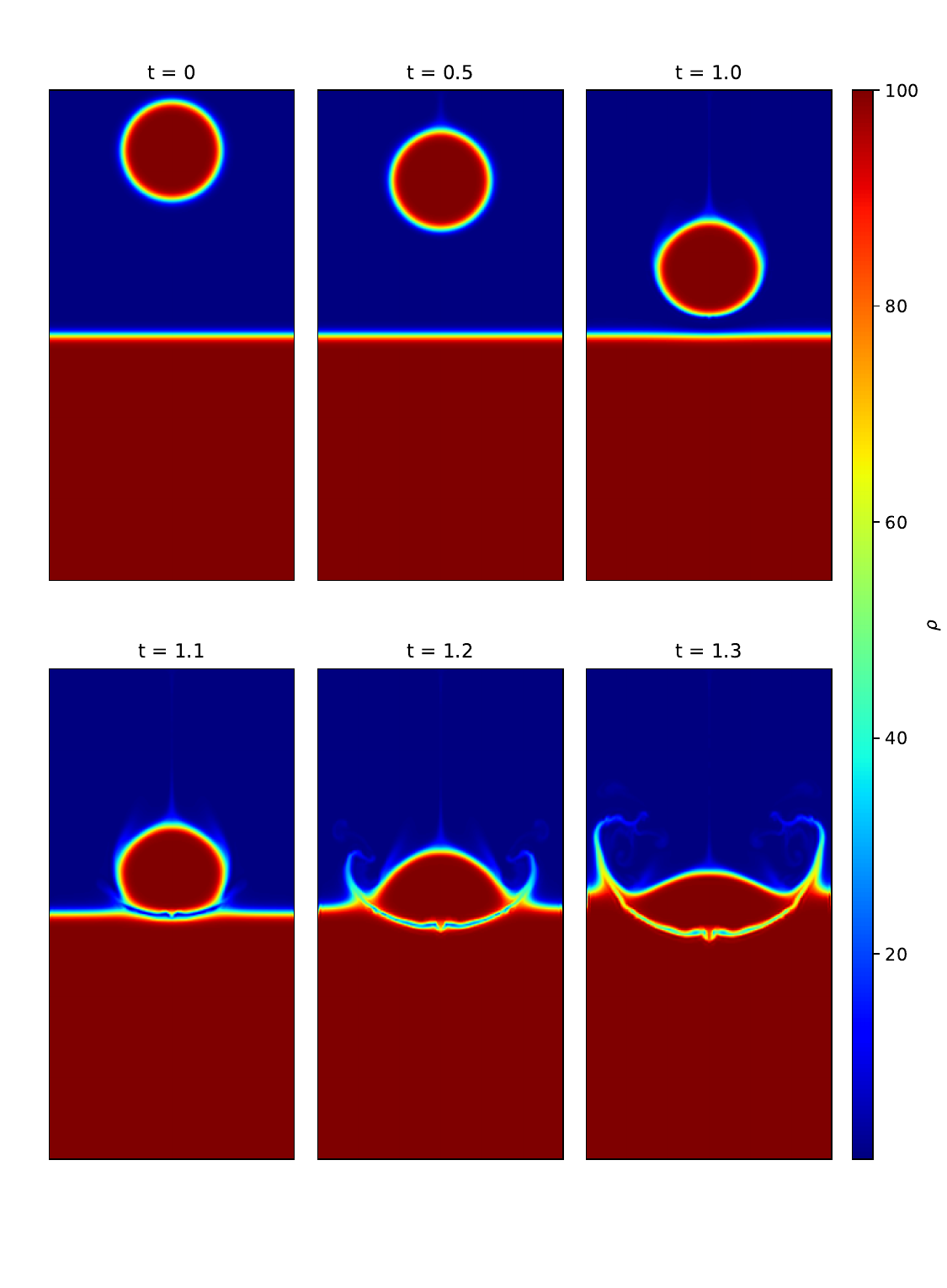}
\caption{Falling-bubble density evolution on the structured \(50\times100\)
RT mesh at different times.}
\label{fig:falling-density-panel}
\end{figure}

\begin{figure}[!htb]
\centering
\begin{minipage}{0.32\textwidth}
\centering
\includegraphics[width=\linewidth]{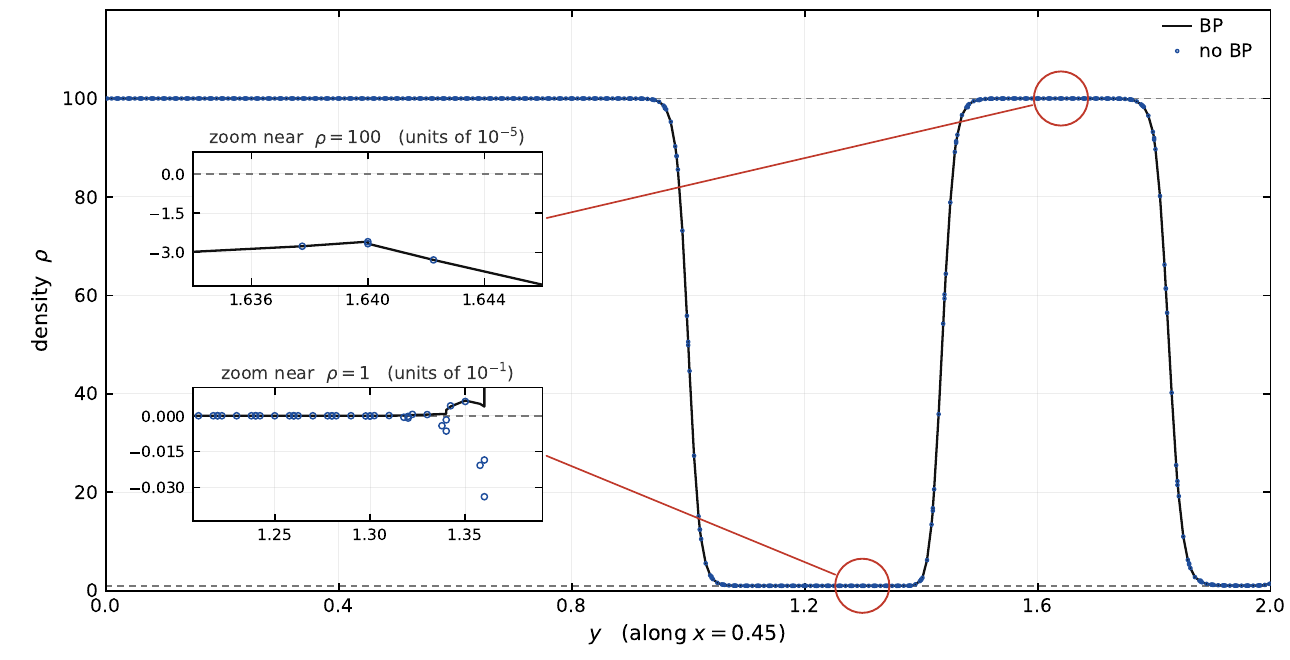}\par
\small (a) \(x=0.45\)
\end{minipage}\hfill
\begin{minipage}{0.32\textwidth}
\centering
\includegraphics[width=\linewidth]{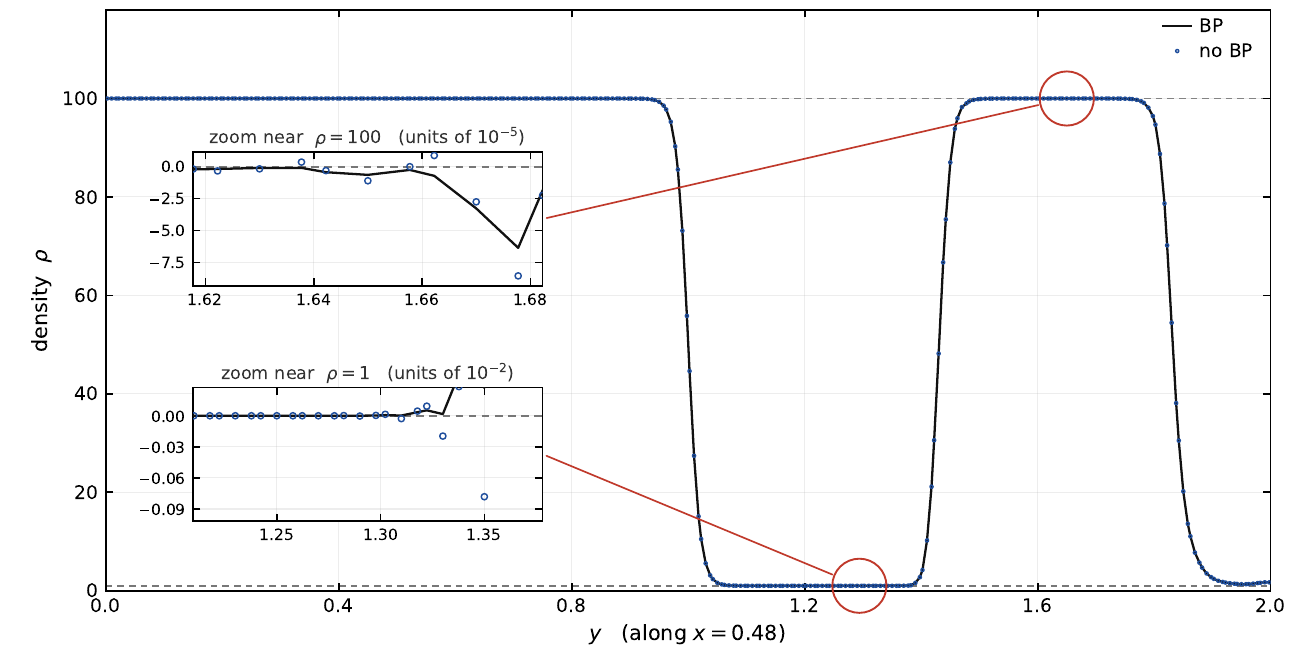}\par
\small (b) \(x=0.48\)
\end{minipage}\hfill
\begin{minipage}{0.32\textwidth}
\centering
\includegraphics[width=\linewidth]{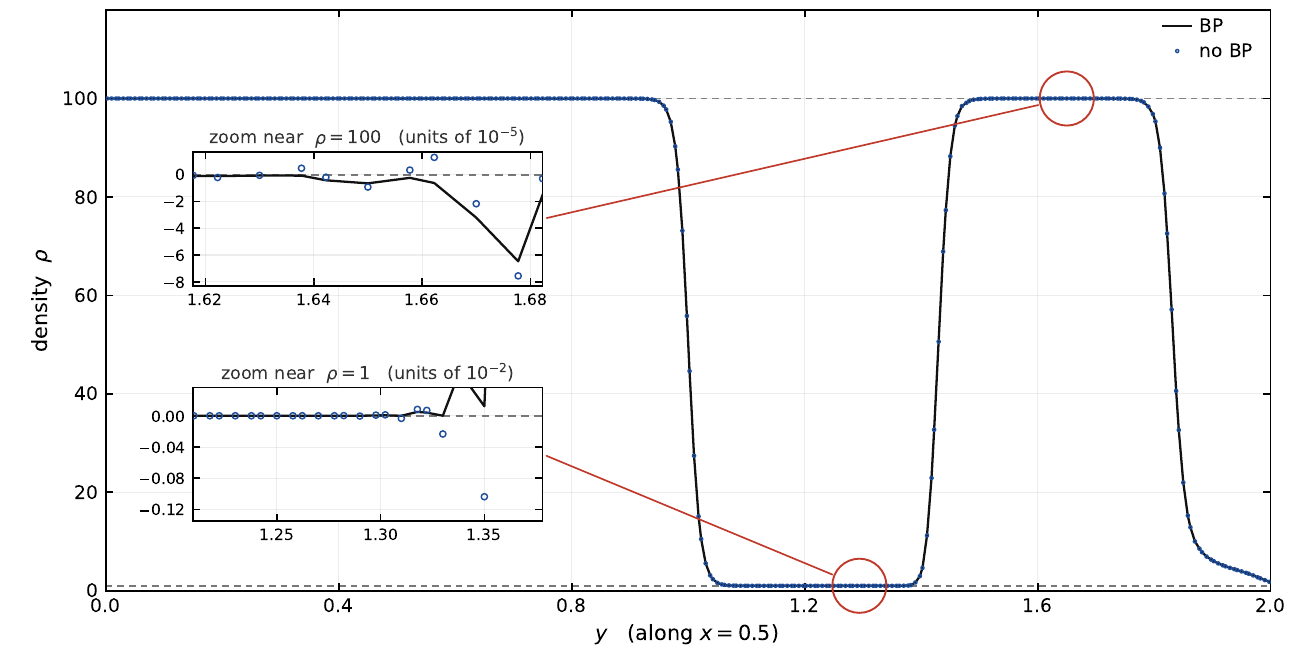}\par
\small (c) \(x=0.50\)
\end{minipage}
\caption{Density profiles for the falling-bubble test along the vertical
lines \(x=0.45,0.48,0.50\), taken at \(t=0.5\), on the structured
\(50\times100\) RT mesh.  The solid lines denote the solutions with BP
limiting and the open circles the solutions without BP limiting; the dashed
lines mark the prescribed density bounds.}
\label{fig:falling-bp-profiles-combined}
\end{figure}

Figure~\ref{fig:falling-diagnostics} reports the divergence, mass-error,
and total-energy histories, and Figure~\ref{fig:falling-density-panel}
shows the descent and deformation of the heavy drop.  
Figure~\ref{fig:falling-bp-profiles-combined} compares the limited and the
unlimited density profiles along the vertical lines \(x=0.45\), \(0.48\),
and \(0.50\).

\FloatBarrier
\subsubsection{All-wall-slip lock exchange}

The lock-exchange geometry and initial data follow
Lundgren and Nazarov~\citep{LundgrenNazarov2024JCP}, with
\(\Omega=[0,32]\times[-0.5,0.5]\), \(\rho_1=1\), \(\rho_2=0.7\), gate
position \(x_d=14\), and \(\vect{u}_0=\vect{0}\).  The regularized initial
density is
\begin{equation}
  \rho_0(x,y)
  =
  \frac{\rho_1+\rho_2}{2}
  -
  \frac{\rho_1-\rho_2}{2}
  \operatorname{erf}\!\left((x-x_d)\sqrt{\mathrm{Re}}\right).
  \label{eq:lock-exchange-initial-density}
\end{equation}
With tank height \(H=1\) and gravitational magnitude \(g_0=1\), the
reference scales are
\begin{equation}
  U_* = \sqrt{\frac{g_0(\rho_1-\rho_2)H}{\rho_1}},
  \qquad
  \mathrm{Re}=\frac{U_*H}{\nu_{\mathrm{kin}}},
  \label{eq:lock-exchange-reference-scales}
\end{equation}
where \(\nu_{\mathrm{kin}}\) denotes the kinematic viscosity.
We set \(\mathrm{Re}=4000\) and \(g=-1\).
For this test, the viscous term in
\eqref{eq:algfull-conservative-system} is replaced by
\(-\nabla\cdot(\mu(\rho)\nabla\vect{u})\), with
\(\mu(\rho)=\rho\nu_{\mathrm{kin}}\).
In each momentum solve, the viscosity is evaluated from
\(\rho_h^{n+1}\), and the gravitational load is
\(\vect{f}^{n+1}=\rho_h^{n+1}(0,g)^{\mathsf T}\).

\begin{figure}[!htb]
\centering
\includegraphics[width=\textwidth]{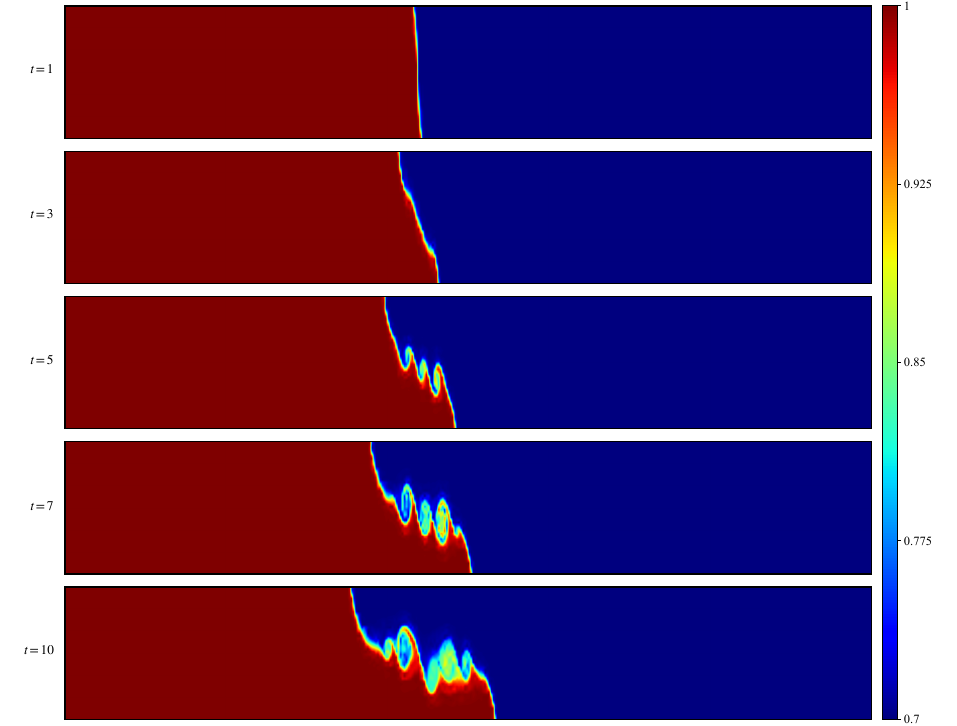}
\caption{All-wall-slip lock-exchange density evolution computed with the
\(\mathrm{RT}_2/\mathbb Q_2\) velocity--pressure pair and discontinuous
\(\mathbb Q_2\) density on the \(N=20\) structured quadrilateral mesh at
\(t=1,3,5,7,10\), ordered from top to bottom.  All panels use the common
display range \(0.7\le\rho_h\le1\).}
\label{fig:lock-exchange-rt-slip-density}
\end{figure}

\begin{figure}[!htb]
\centering
\includegraphics[width=\textwidth]{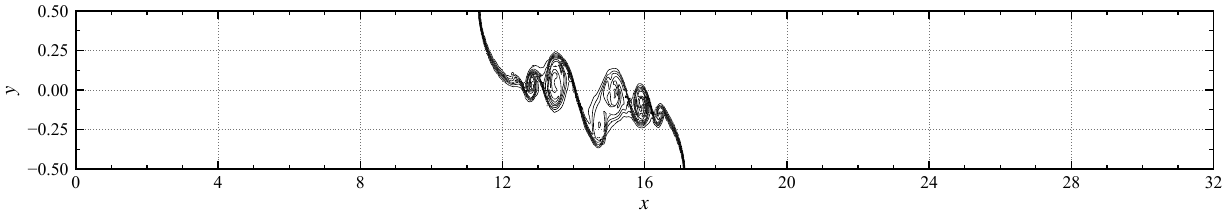}
\caption{Final density contours for the all-wall-slip RT lock exchange on the
\(N=20\) mesh at \(t=10\).}
\label{fig:lock-exchange-rt-slip-contours}
\end{figure}

The discretization uses \(\mathrm{RT}_2/\mathbb Q_2\) for velocity and
pressure and discontinuous \(\mathbb Q_2\) for density on a structured
\(N=20\) mesh
(\(32N\times N\) rectangles of side length \(1/N\)), with
\(\Delta t=1/110\).  All walls are impermeable and satisfy slip boundary
conditions: \(\vect{u}_h\cdot\vect{n}=0\), with the natural tangential
viscous condition (\(\Gamma_D=\varnothing\)).
Figures~\ref{fig:lock-exchange-rt-slip-density} and
\ref{fig:lock-exchange-rt-slip-contours} show the two gravity currents and
the interfacial roll-up; the comparison is qualitative
\citep{LundgrenNazarov2024JCP}.

\begin{figure}[!htb]
\centering
\begin{minipage}{0.32\textwidth}
\centering
\includegraphics[width=\textwidth]{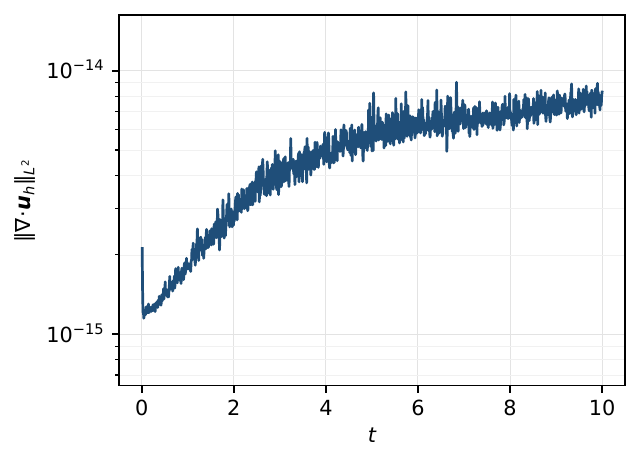}\\
\small (a) divergence
\end{minipage}\hfill
\begin{minipage}{0.32\textwidth}
\centering
\includegraphics[width=\textwidth]{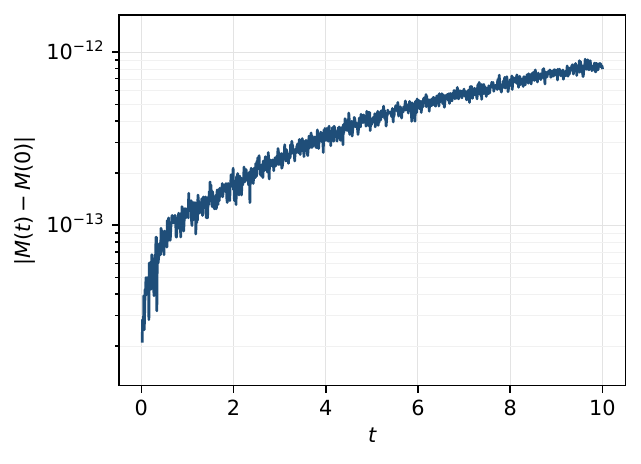}\\
\small (b) absolute mass error
\end{minipage}\hfill
\begin{minipage}{0.32\textwidth}
\centering
\includegraphics[width=\textwidth]{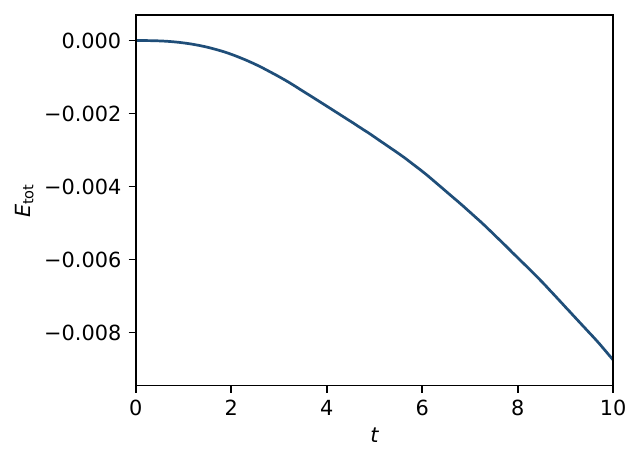}\\
\small (c) total energy
\end{minipage}
\caption{Diagnostics for the all-wall-slip RT lock-exchange calculation:
(a) \(\|\nabla\cdot \vect{u}_h\|_{L^2}\) and (b) the absolute mass deviation
\(|\mathcal M_h(t)-\mathcal M_h(0)|\), both on base-10 logarithmic vertical axes, and (c) total
energy with potential energy measured relative to \(y=0\).}
\label{fig:lock-exchange-rt-slip-diagnostics}
\end{figure}

Figure~\ref{fig:lock-exchange-rt-slip-diagnostics} reports the divergence,
mass error, and total energy defined in \eqref{eq:rt-energy-definition},
with potential energy measured relative to \(y=0\).

\subsection{Discussion}
\label{sec:discussion}

The numerical results support the method's design.  The mixed BDM/RT
discretization supplies transport velocities with the local incompressibility
and normal-trace properties required by the density analysis, while the
partitioned SSPRK--BDF integrator separates the explicit bound-preserving
transport update from the linear mixed momentum solve.  This coupling is
particularly useful when the density range is wide: the limiter controls the
density polynomial without altering its cell average, and the momentum solve
continues to use the updated conservative density.

The guarantees should be interpreted within their stated scope.  The
rigorous bound result is for source-free density transport with periodic or
impermeable boundaries, admissible starting histories, the specified
quadrature and polynomial-degree assumptions, and a CFL restriction applied
to each extrapolated stage velocity.  The limiter controls the cell averages
and the selected control-point values; it does not by itself establish a
pointwise bound at every point of an element.  The rectangular RT results
use the corresponding tensor-product construction and are supported here by
the numerical verification.

The partitioned treatment also introduces practical considerations.  The
stagewise CFL restriction must account for extrapolated velocities, and the
startup procedure must produce both admissible density states and solenoidal
velocity histories.  In addition, although the total-energy histories in
the benchmark calculations are useful diagnostics, a fully discrete energy
estimate for the coupled BP-limited method has not been established.  These
issues identify natural directions for further analysis, including energy
stability, sharper timestep conditions, and extensions to more general
constitutive laws and three-dimensional geometries.

\section{Conclusion}
\label{sec:conclusion}

We have developed a structure-preserving discontinuous Galerkin method for
variable-density incompressible flow with a partitioned SSPRK--BDF time
integrator.  The BDM/RT mixed momentum discretization and velocity
extrapolation provide exactly divergence-free transport velocities with
single-valued normal traces, while conservative DG transport and
mean-preserving limiting enforce density bounds and global mass conservation
under the stated assumptions.
The analysis establishes the cell-average and control-point bounds for
source-free transport on triangular meshes, and the numerical results verify
the expected convergence, divergence-free behavior, mass conservation, and
bound preservation on both triangular and structured rectangular meshes.
The method remains effective in low-density and large-density-ratio tests.

\bibliographystyle{abbrv}
\bibliography{references}

\end{document}